\documentclass[11pt]{amsart}

\usepackage{amsmath} 
\usepackage{amssymb}

\newcommand{\V}{{\bf V}}

\newcommand{\cf}{{\rm cf}}

\newcommand{\bB}{{\mathbb B}}
\newcommand{\bQ}{{\mathbb Q}}
\newcommand{\bP}{{\mathbb P}}

\newcommand{\cA}{{\mathcal A}}
\newcommand{\cD}{{\mathcal D}}
\newcommand{\rest}{\restriction}
\newcommand{\otp}{{\rm otp}}
\newcommand{\cl}{{\rm cl}}
\newcommand{\qlm}{{{\mathbb Q}^{\mu,\lambda}_n}}
\newcommand{\bPk}{{\mathbb P}_\kappa}
\newcommand{\bPz}{{\mathbb P}^0_\kappa}
\newcommand{\bPj}{{\mathbb P}^1_\kappa}

\newcommand{\qs}{{\mathbb Q}_S}
\newcommand{\comp}{\circ}
\newcommand{\forces}{\Vdash}
\newcommand{\ut}{{\rm ut}}
\newcommand{\tig}{{\rm t}}
\newcommand{\ind}{{\rm ind}}
\newcommand{\irr}{{\rm irr}}
\newcommand{\hd}{{\rm hd}}
\newcommand{\hL}{{\rm hL}}
\newcommand{\inc}{{\rm inc}}
\newcommand{\Id}{{\rm Id}}
\newcommand{\Sub}{{\rm Sub}}
\newcommand{\aut}{{\rm Aut}}
\newcommand{\Dep}{{\rm Depth}}
\newcommand{\inv}{{\rm inv}}

\newcommand{\Ult}{{\rm Ult}}

\newcommand{\cK}{{\mathcal K}}
\newcommand{\xs}{{\mathcal X}_S}

\newcommand{\QED}{\hfill\vrule width 6pt height 6pt depth 0pt\vspace{0.1in}}
\newcommand{\Proof}{\noindent{\sc Proof} \hspace{0.2in}} 

\newtheorem{theorem}{Theorem}[section] 
\newtheorem{claim}{Claim}[theorem]
\newtheorem{lemma}[theorem]{Lemma} 
\newtheorem{proposition}[theorem]{Proposition}

\theoremstyle{definition}
\newtheorem{definition}[theorem]{Definition}
\newtheorem{problem}[theorem]{Problem}

\theoremstyle{remark}
\newtheorem{conclusion}[theorem]{Conclusion}
\newtheorem{remark}[theorem]{Remark}
\newtheorem{notation}[theorem]{Notation}

\title{More on cardinal invariants of Boolean algebras}

\author{Andrzej Ros{\l}anowski}
\address{Institute of Mathematics\\
The Hebrew University\\
91 904 Jerusalem, Israel\\
and Mathematical Institute\\
Wroc{\l}aw University\\
50 384 Wroc{\l}aw, Poland}
\email{roslanow@math.idbsu.edu}
\urladdr{http://math.idbsu.edu/$\sim$roslanow}

\author{Saharon Shelah}
\address{Institute of Mathematics\\
 The Hebrew University of Jerusalem\\
 91904 Jerusalem, Israel\\
 and  Department of Mathematics\\
 Rutgers University\\
 New Brunswick, NJ 08854, USA}
\email{shelah@math.huji.ac.il}
\urladdr{http://www.math.rutgers.edu/$\sim$shelah}
\thanks{The research of the second author was partially supported by The
Israel Science Foundation. Publication 599} 

\begin{document}

\begin{abstract}
We address several questions of Donald Monk related to irredundance and spread
of Boolean algebras, gaining both some ZFC knowledge and consistency results.
We show in ZFC that $\irr(\bB_0\times\bB_1)=\max\{\irr(\bB_0),\irr(\bB_1)\}$.
We prove consistency of the statement ``there is a Boolean algebra $\bB$ such
that $\irr(\bB)<s(\bB\circledast\bB)$'' and we force a superatomic Boolean
algebra $\bB_*$ such that $s(\bB_*)=\inc(\bB_*)=\kappa$, $\irr(\bB_*)=\Id(
\bB_*)=\kappa^+$ and $\Sub(\bB_*)=2^{\kappa^+}$. Next we force a superatomic
algebra $\bB_0$ such that $\irr(\bB_0)<\inc(\bB_0)$ and a superatomic algebra
$\bB_1$ such that $\tig(\bB_1)> {\rm Aut}(\bB_1)$. Finally we show that
consistently there is a Boolean algebra $\bB$ of size $\lambda$ such that
there is no free sequence in $\bB$ of length $\lambda$, there is an
ultrafilter of tightness $\lambda$ (so $\tig(\bB)=\lambda$) and $\lambda\notin
\Dep_{\rm Hs}(\bB)$.
\end{abstract}

\maketitle

\section{Introduction}
In the present paper we answer (sometimes partially only) several questions of
Donald Monk concerning cardinal invariants of Boolean algebras. Most of our
results are consistency statements, but we get some ZFC knowledge too.

For a systematic study and presentation of current research on cardinal
invariants of Boolean algebras (as well as for a long list of open problems)
we refer the reader to Monk \cite{M2}. Some of the relevant definitions are
listed at the end of this section.
\medskip

\noindent{\bf Content of the paper:}\qquad In the first section we show that
the difference between $s_n(\bB)$ and $s_N(\bB)$ (for $n<N$) can be reasonably
large, with the only restriction coming from the inequality $s_n(\bB)\geq
2^{s_N}(\bB)$ (a consistency result; for the definitions of the invariants
see below). It is relevant for the description of the behaviour of spread in
ultraproducts: we may conclude that it is consistent that $s(\prod\limits_{n
\in\omega}\bB_n/D)$ is much larger than $\prod\limits_{n\in\omega}s(\bB_n)/D$. 
In the following section we answer \cite[Problem 24]{M2} showing that $\irr(
\bB_0\times\bB_1)=\max\{\irr(\bB_0),\irr(\bB_1)\}$ (a ZFC result). A partial
answer to \cite[Problem 27]{M2} is given in the third section, where we show
that, consistently, there is a Boolean algebra $\bB$ such that $\irr(\bB)<s(\bB
\circledast\bB)$. In particular, this shows that the parallel statement to the
result of section 2 for free product may fail. Note that proving the result of
section 3 in ZFC is a really difficult task, as so far we even do not know if
(in ZFC) there are Boolean algebras $\bB$ satisfying $\irr(\bB)<|\bB|$. In
section 4 we force a superatomic Boolean algebra $\bB$ such that $s(\bB)=\inc(
\bB)=\kappa$, $\irr(\bB)=\Id(\bB)=\kappa^+$ and $\Sub(\bB)=2^{\kappa^+}$. This
gives answers to \cite[Problems 73, 77, 78]{M2} as stated (though the problems
in ZFC remain open). Next we present some modifications of this forcing notion
and in the fifth section we answer \cite[Problems 79, 81]{M2} forcing
superatomic Boolean algebras $\bB_0,\bB_1$ such that $\irr(\bB_0)<\inc(\bB_0)$
and $\aut(\bB_1)<\tig(\bB_1)$. Finally in the last section we show that
(consistently) there is a Boolean algebra $\bB$ of size $\lambda$ such that
there is no free sequence in $\bB$ of length $\lambda$, there is an
ultrafilter in $\Ult(\bB)$ of tightness $\lambda$ (so $\tig(\bB)=\lambda$) and
$\lambda\notin\Dep_{\rm Hs}(\bB)$. This gives answers to \cite[Problems 13,
41]{M2}. Lastly we use one of the results of \cite{Sh:233} to show that
$2^{\cf(\tig(\bB))}<\tig(\bB)$ implies $\tig(\bB)\in\Dep_{\rm Hs}(\bB)$.
\medskip

\noindent{\bf Notation:}\qquad Our notation is rather standard and
compatible with that of classical textbooks on set theory (like Jech
\cite{J}) and Boolean algebras (like Monk \cite{M1}, \cite{M2}). However
in forcing considerations we keep the older tradition that
\begin{center}
{\em
the stronger condition is the greater one
}
\end{center}
Let us list some of our notation and conventions.

\begin{notation}
\begin{enumerate}
\item A name for an object in a forcing extension is denoted with a dot above
(like $\dot{X}$) with one exception: the canonical name for a generic filter
in a forcing notion $\bP$ will be called $\Gamma_\bP$.
\item $\alpha,\beta,\gamma,\delta,\ldots$ will denote ordinals and
$\kappa,\mu,\lambda,\theta$ will stand for (always infinite) cardinals.
\item For a set $X$ and a cardinal $\lambda$, $[X]^{\textstyle<\lambda}$
stands for the family of all subsets of $X$ of size less than $\lambda$. If
$X$ is a set of ordinals then its order type is denoted by $\otp(X)$.
\item In Boolean algebras we use $\vee$ (and $\bigvee$), $\wedge$ (and
$\bigwedge$) and $-$ for the Boolean operations. If $\bB$ is a Boolean
algebra, $x\in\bB$ then $x^0=x$, $x^1=-x$. The Stone space of the algebra
$\bB$ is called $\Ult(\bB)$.
\item For a subset $Y$ of an algebra $\bB$, the subalgebra of $\bB$ generated
by $Y$ is denoted by $\langle Y\rangle_{\bB}$.
\item The sign $\circledast$ stands for the operation of the free product of
Boolean algebras and the product is denoted by $\times$.
\end{enumerate}
\end{notation}
\smallskip

\noindent{\bf The invariants:}\qquad Below we recall some definitions and
formalism from \cite{RoSh:534} (see \cite{M2} too).

\begin{definition}
For a (not necessary first order) theory $T$ in the language of Boolean
algebras plus one distinguished unary predicate $P_0$ plus, possibly, some
others $P_1,P_2,\ldots$ we define cardinal invariants $\inv_T$, $\inv^+_T$ of
Boolean algebras by (for a Boolean algebra $\bB$):
\begin{quotation}
$\inv_T(\bB)\stackrel{\rm def}{=}\sup\{|P_0|: (B,P_n)_n \mbox{ is a model of }
T\}$, 

$\inv^+_T(\bB)\stackrel{\rm def}{=}\sup\{|P_0|^+: (B,P_n)_n\mbox{ is a model
of } T\}$. 
\end{quotation}
\end{definition}

We think of the spread $s(\bB)$ of a Boolean algebra $\bB$ as 
\begin{enumerate}
\item[$(\otimes_s)$] $s(\bB)=\sup\{|X|: X\subseteq B\mbox{ is
ideal-independent}\}$   
\end{enumerate}
(it is one of the equivalent definitions, see \cite[Thm 13.1]{M2}). Thus we
can write $s(\bB)=s_\omega(\bB)$, where 
\begin{definition}
\begin{enumerate}
\item $\phi_n^s$ is the formula saying that no member of $P_0$ can be
covered by union of $n+1$ other elements of $P_0$.
\item For $0<n\leq\omega$ let $T^n_s=\{\phi_k^s:k<n\}$.
\item For a Boolean algebra $\bB$ and $0<n\leq\omega$:\quad $s^{(+)}_n(\bB)=
\inv_{T^n_s}^{(+)}(\bB)$.
\end{enumerate}
\end{definition}
The hereditary density $\hd(\bB)$ and the hereditary Lindel\"of degree
$\hL(\bB)$ of a Boolean algebra $\bB$ are treated in a similar manner. We use 
\cite[Thm 16.1]{M2} and \cite[Thm 15.1]{M2} to define them as
\begin{enumerate}
\item[$(\otimes_{\hd})$] $\hd(\bB)=\sup\{|\kappa|:$ there is a strictly
decreasing sequence of ideals (in $\bB$) of the length $\kappa\ \}$,
\item[$(\otimes_{\hL})$] $\hL(\bB)=\sup\{|\kappa|:$ there is a strictly
increasing sequence of ideals (in $\bB$) of the length $\kappa\ \}$.
\end{enumerate}
This leads us directly to the following definition.
\begin{definition}
\begin{enumerate}
\item[1.\ ] Let the formula $\psi$ say that $P_1$ is a well ordering of
$P_0$ (denoted by $<_1$).
\item[2.\ ] For $n<\omega$ let $\phi^{\hd}_n$, $\phi^{\hL}_n$ be the
following formulas:
\item[$\phi^{\hd}_n\equiv$] $\psi\ \&\ (\forall x_0,\ldots,x_{n+1}\in P_0)(
x_0<_1\ldots<_1 x_{n+1}\ \Rightarrow\ x_0\not\leq x_1\vee\ldots\vee x_{n+1})$
\item[$\phi^{\hL}_n\equiv$] $\psi\ \&\ (\forall x_0,\ldots,x_{n+1}\in P_0)(
x_{n+1}<_1\ldots<_1 x_0\ \Rightarrow\ x_0\not\leq x_1\vee\ldots\vee x_{n+1})$.
\item[3.\ ] For $0<n\leq\omega$ we let $T^n_{\hd}=\{\phi^{\hd}_k:k<n\}$,
$T^n_{\hL}=\{\phi^{\hL}_k: k<n\}$. 
\item[4.\ ] For a Boolean algebra $\bB$ and $0<n\leq\omega$:
\[\hd^{(+)}_n(\bB)=\inv^{(+)}_{T^n_{\hd}}(\bB),\quad\quad \hL^{(+)}_n(\bB)=
\inv^{(+)}_{T^n_{\hL}}(\bB).\]
\end{enumerate}
\end{definition}

We use the following characterization of tightness (see \cite[\S 12]{M2}): 
\[\tig(\bB)=\sup\{|\alpha|:\mbox{ there exists a free sequence of the length
}\alpha\mbox{ in }\bB\}.\] 

\begin{definition}
\begin{enumerate}
\item Let $\psi$ be the sentence saying that $P_1$ is a well ordering of $P_0$
(we denote the respective order by $<_1$). For $k,l<\omega$ let
$\phi_{k,l}^\tig$ be the sentence asserting that 
\begin{quotation}
\noindent for each $x_0,\ldots,x_k,y_0,\ldots,y_l\in P_0$

\noindent if $x_0<_1\ldots<_1 x_k<_1 y_0<_1\ldots <_1 y_l$ then
$\bigwedge\limits_{i\leq k}x_i\not\leq\bigvee\limits_{i\leq l} y_i$,
\end{quotation}
and let the sentence $\phi_{k,l}^{\ut}$ say that
\[\mbox{for each distinct }x_0,\ldots,x_k,y_0,\ldots,y_l\in P_0\mbox{ we have }
\bigwedge\limits_{i\leq k}x_i\not\leq\bigvee\limits_{i\leq l} y_i.\]
\item For $n,m\leq\omega$ let $T^{n,m}_\tig=\{\phi_{k,l}^\tig:k<n,l<m\}\cup
\{\psi\}$ and $T^{n,m}_{\ut}=\{\phi^{\ut}_{k,l}: k<n,l<m\}$ and for a Boolean
algebra $\bB$: 
\[\tig_{n,m}(\bB)=\inv_{T^{n,m}_\tig}(\bB)\qquad\&\qquad\ut_{n,m}(\bB)=
\inv_{T^{n,m}_{\ut}}(\bB).\]
\end{enumerate}
\end{definition}

The irredundance $\irr(\bB)$ of a Boolean algebra $\bB$ is the supremum of
cardinalities of sets $X\subseteq\bB$ such that $(\forall x\in X)(x\notin
\langle X\setminus\{x\}\rangle_{\bB})$. 
\begin{definition}[compare {\cite[p. 144]{M2}}]
Let $n\leq\omega$ and let $T^n_{\irr}$ be the theory of the language of
Boolean algebras plus a predicate $P_0$, which says that for each $m<n$ and
a Boolean term $\tau(y_0,\ldots,y_m)$ we have
\[(\forall x\in P_0)(\forall x_0,\ldots,x_m\in P_0\setminus\{x\})(x\neq \tau
(x_0,\ldots,x_m)).\]
For Boolean algebra $\bB$ we define $\irr^{(+)}_n(\bB)=\inv^{(+)}_{T^n_{
\irr}}(\bB)$ (so $\irr_\omega(\bB)=\irr(\bB)$).
\end{definition}
The incomparability number $\inc(\bB)$ is the supremum of cardinalities of
sets of pairwise incomparable elements. The number of ideals in $\bB$ is
denoted by $\Id(\bB)$, $\aut(\bB)$ stands for the number of automorphisms of
the algebra $\bB$, and the number of subalgebras of $\bB$ is denoted by $\Sub(
\bB)$.  

\section{Forcing for spread}
The aim of this section is to show that for $N$ much larger than $n$, the
inequalities $2^{s_N}(\bB)\geq s_n(\bB)\geq s_N(\bB)$ (see \cite[Thm
13.6]{M2}) seem to be the only restriction on the jumps between $s_N$ and
$s_n$. The forcing notion defined in \ref{1.1}(2) below is a modification of
the one from \cite[\S 2]{Sh:479} and a relative of the forcing notion from
\cite[\S 15]{Sh:620}.

\begin{definition}
\label{1.1}
1)\quad For a set $w$ and a family $F\subseteq 2^{\textstyle w}$ we define

\noindent $\cl(F)=\{g\in 2^{\textstyle w}: (\forall u\in [w]^{\textstyle
<\omega})(\exists f\in F)(f\rest u=g\rest u)$,

\noindent $\bB_{(w,F)}$ is the Boolean algebra generated freely by
$\{x_\alpha:\alpha\in w\}$ except that

if $u_0,u_1\in [w]^{\textstyle <\omega}$ and there is no $f\in F$ such
that $f\rest u_0\equiv 0$, $f\rest u_1\equiv 1$

then $\bigwedge\limits_{\alpha\in u_1} x_\alpha\wedge \bigwedge\limits_{
\alpha\in u_0} (-x_\alpha)=0$.

\noindent 2)\quad Let $\mu\leq\lambda$ be cardinals, $0<n<\omega$. We define
forcing notion $\qlm$: 
\smallskip

\noindent {\bf a condition}\quad is a pair $p=(w^p,F^p)$ such that
$w^p\in [\lambda]^{\textstyle <\mu}$, $F^p\subseteq 2^{\textstyle w^p}$,
$|F^p|<\mu$ and for every $u\in [w^p]^{\textstyle \leq\! n}$ there is $f^*:
w^p\setminus u\longrightarrow 2$ such that

\qquad if $h:u\longrightarrow 2$ then $f^*\cup h\in F^p$;

\noindent {\bf the order}\quad is given by\quad $p\leq q$ if and only if
$w^p\subseteq w^q$ and
\[(\forall f\in F^q)(f\rest w^p\in \cl(F^p))\ \mbox{ and }\ (\forall f\in F^p)
(\exists g\in F^q)(f\subseteq g).\]
\end{definition}

\begin{proposition}
[see {\cite[2.6]{Sh:479}}]
\label{1.2}
\begin{enumerate}
\item If $p\in \qlm$, $f\in F^p$ then $f$ extends to a homomorphism from
$\bB_p$ to $\{0,1\}$ (i.e.~it preserves the equalities from the
definition of $\bB_p$).
\item If $p\in\qlm$, $\tau(y_0,\ldots,y_k)$ is a Boolean term and
$\alpha_0,\ldots,\alpha_k\in w^p$ are distinct then

$\bB_p\models\tau(x_{\alpha_0},\ldots,x_{\alpha_k})\neq 0$\quad if and
only if 
\[(\exists f\in F^p)(\{0,1\}\models\tau(f(\alpha_0),\ldots,f(\alpha_k))=1).\]
\item If $p,q\in\qlm$, $p\leq q$ then $\bB_p$ is a subalgebra of $\bB_q$.\QED
\end{enumerate}
\end{proposition}

\begin{proposition}
\label{1.3}
Assume $\mu^{<\mu}=\mu\leq\lambda$, $0<n<\omega$. Then
\begin{enumerate}
\item $\qlm$ is a $\mu$--complete forcing notion of size $\lambda^{<\mu}$,
\item $\qlm$ satisfies $\mu^+$-cc.
\end{enumerate}
\end{proposition}

\Proof This is almost exactly like \cite[2.7]{Sh:479}. For (1) no changes are
required; for (2) one has to check that the condition defined as there is
really in $\qlm$. So suppose that $\langle p_\alpha:\alpha<\mu^+\rangle
\subseteq\qlm$. Applying standard ``cleaning procedure'' find $\alpha_0<
\alpha_1<\mu^+$ such that
\begin{itemize}
\item $\otp(w^{p_{\alpha_0}})=\otp(w^{p_{\alpha_1}})$,
\item if $H:w^{p_{\alpha_0}}\longrightarrow w^{p_{\alpha_1}}$ is the order
preserving mapping then $H\rest (w^{p_{\alpha_0}}\cap w^{p_{\alpha_1}})$ is
the identity on $w^{p_{\alpha_0}}\cap w^{p_{\alpha_1}}$ and $F^{p_{\alpha_0}}
=\{f\comp H: f\in F^{p_{\alpha_1}}\}$
\end{itemize}
(remember $\mu^{<\mu}=\mu$; use $\Delta$--lemma). Let $w^q=w^{p_{\alpha_0}}
\cup w^{p_{\alpha_1}}$ and
\[F^q=\{f\cup g: f\in F^{p_{\alpha_0}}\ \&\ g\in F^{p_{\alpha_1}}\ \&\
f\rest (w^{p_{\alpha_0}}\cap w^{p_{\alpha_1}}) = g\rest (w^{p_{\alpha_0}}
\cap w^{p_{\alpha_1}})\}.\]
To check that $q=(w^q,F^q)$ is in $\qlm$ suppose that $u\in [w^q]^{\textstyle
\leq\! n}$ and let $u^*=H^{-1}[u\cap w^{p_{\alpha_1}}]\cup (u\cap w^{p_{
\alpha_0}})\in [w^{p_{\alpha_0}}]^{\textstyle \leq\! n}$. Let $f^*_0:
w^{p_{\alpha_0}}\setminus u^*\longrightarrow 2$ be such that\quad if $h:u^*
\longrightarrow 2$ then $f^*_0\cup h\in F^{p_{\alpha_0}}$. Next, let $f^*:
w^{p_{\alpha_0}}\setminus u\longrightarrow 2$ be such that $f^*_0\subseteq
f^*$ and if $\alpha\in u^*\setminus u$ then $f^*(\alpha)=0$, and let $g^*:
w^{p_{\alpha_1}}\setminus u\longrightarrow 2$ be such that $f^*_0\comp H^{-1}
\subseteq g^*$ and if $\alpha\in H[u^*]\setminus u$ then $g^*(\alpha)=0$. Now
it should be clear that 
\[\mbox{if }\quad h:u\longrightarrow 2\quad\mbox{ then }\quad (f^*\cup g^*)
\cup h\in F^q.\]
Verifying that both $p_{\alpha_0}\leq q$ and $p_{\alpha_1}\leq q$ is
even easier. \QED
\medskip

Let $\dot{\bB}$ be the $\qlm$--name for $\bigcup\{\bB_p:p\in\Gamma_{\qlm}
\}$. It follows from \ref{1.2} that  
\[\forces_{\qlm}\mbox{`` $\dot{\bB}$ is a Boolean algebra generated by
$\{x_\alpha:\alpha<\lambda\}$ ''}\]
and, for a condition $p\in\qlm$,
\[p\forces_{\qlm}\mbox{`` }\langle x_\alpha:\alpha\in w^p\rangle_{\dot{
\bB}}=\bB_p\mbox{ ''}.\]

\begin{theorem}
\label{1.4}
Assume $\mu^{<\mu}=\mu\leq\lambda$ and $0<N,n<\omega$ are such that $2^{n/2}+
n\leq N$. Then
\[\forces_{\qlm}\mbox{`` }\ind^+_n(\dot{\bB})=\lambda^+\ \mbox{ and }\
\tig^+_{1,N}(\dot{\bB})=\tig^+_{N,1}(\dot{\bB})=\ind^+(\dot{\bB})=\mu^+
\mbox{ ''.}\] 
\end{theorem}

\Proof It follows immediately from the definition of $\qlm$ (by density
arguments, remembering \ref{1.2}) that
\[\forces_{\qlm}\mbox{`` the sequence $\langle x_\alpha:\alpha<\lambda
\rangle$ is $n$--independent ''.}\]
Suppose now that $\langle \dot{a}_\beta:\beta<\mu^+\rangle$ is a $\qlm$--name
for a $\mu^+$--sequence of elements of $\dot{\bB}$, $p\in\qlm$. For each
$\beta<\mu^+$ choose a condition $p_\beta\geq p$, a Boolean term $\tau_\beta$
and ordinals $\bar{\alpha}(\beta,0)<\ldots<\bar{\alpha}(\beta,\ell_\beta)<
\lambda$ such that
\[p_\beta\forces_{\qlm}\dot{a}_\beta=\tau_\beta(x_{\bar{\alpha}(\beta,0)},
\ldots,x_{\bar{\alpha}(\beta,\ell_\beta)}).\]
By $\Delta$--system arguments, passing to a subsequence and increasing
$p_\beta$'s, we may assume that
\begin{enumerate}
\item[(i)]\ \ \ $\tau_\beta=\tau$, $\ell_\beta=\ell$ and $\bar{\alpha}(\beta,
0),\ldots,\bar{\alpha}(\beta,\ell)\in w^{p_\beta}$,
\item[(ii)]\ \  $\otp(w^{p_{\beta_0}})=\otp(w^{p_{\beta_1}})$ and $\otp(w^{
p_{\beta_0}}\cap\bar{\alpha}(\beta_0,j))=\otp(w^{p_{\beta_1}}\cap\bar{\alpha}
(\beta_1,j))$ for $j\leq \ell$, $\beta_0,\beta_1<\mu^+$,
\item[(iii)]\   $\{w^{p_\beta}:\beta<\mu^+\}$ forms a $\Delta$--system of sets
with heart $w^*$, 
\item[(iv)]\ \  if $H_{\beta_0,\beta_1}:w^{p_{\beta_0}}\longrightarrow w^{p_{
\beta_1}}$ is the order preserving mapping then $H_{\beta_0,\beta_1}\rest w^*$
is the identity on $w^*$ and $F^{p_{\beta_0}}=\{f\comp H_{\beta_0,\beta_1}:
f\in F^{p_{\beta_1}}\}$. 
\end{enumerate}
After this ``cleaning procedure'' look at the conditions $p_0,\ldots,p_N$. We
want to show that they have a common upper bound $q\in\qlm$ such that
$q\forces_{\qlm}\mbox{`` }\dot{a}_0\wedge\bigwedge\limits_{j<N} (-\dot{a}_{1+
j})=0\mbox{ ''}$. To this end define:
\[w^q=w^{p_0}\cup\ldots\cup w^{p_N}\quad\mbox{ and}\]
\[\begin{array}{ll}
F^q=\big\{f_0\cup\ldots\cup f_N:& f_0\in F^{p_0},\ldots,f_N\in F^{p_N},
\ f_0\rest w^*=\ldots=f_N\rest w^*,\\
\ &\mbox{and if }\ \ \{0,1\}\models\tau(f_0(\bar{\alpha}(0,0)),\ldots,f_0(
\bar{\alpha}(0,\ell)))=1\\
\ &\mbox{then for some }j\in[1,N]\\
\ &\{0,1\}\models\tau(f_j(\bar{\alpha}(j,0)),\ldots,f_j(\bar{\alpha}(j,\ell)))
=1\ \big\}. 
\end{array}\]
Let us check that $q=(w^q,F^q)$ is in $\qlm$. Clearly each $f\in F^q$ is a
function from $w^q$ to 2 and $|F^q|<\mu$. Suppose now that $u\in [w^q]^{
\textstyle\leq\! n}$. Let $u^*=u\cap w^*$ and $u^+=\bigcup\limits_{i\leq N} 
H_{i,0}[u\cap w^{p_i}]\in [w^{p_0}]^{\textstyle\leq\! n}$. One of the sets
$u^*$, $u^+\setminus u^*$ has size at most $n/2$, and first we deal with the
case $|u^*|\leq n/2$. Choose $f^*:w^{p_0}\setminus u^+\longrightarrow 2$ such
that $(\forall h:u^+\longrightarrow 2)(f^*\cup h\in F^{p_0})$. For each $v
\subseteq u^*$ choose $h_v:u^+\longrightarrow 2$ such that $h_v\rest v\equiv
1$, $h_v\rest (u^* \setminus v)\equiv 0$ and
\begin{quotation}
\noindent if there is $h:u^+\longrightarrow 2$ satisfying the above demands
and such that\quad $\{0,1\}\models \tau((f^*\cup h)(\bar{\alpha}(0,0)),\ldots,
(f^*\cup h)(\bar{\alpha}(0,\ell)))=1$\\
then $h_v$ has this property.
\end{quotation}
Since $2^{|u^*|}+n\leq N$ we may choose distinct $i_v\in [1,N]$ for $v
\subseteq u^*$ such that $w^{p_{i_v}}\cap u=u^*$. Now we define functions
$f^*_i: w^{p_i}\setminus u\longrightarrow 2$ (for $i\leq N$) as follows:
\begin{itemize}
\item if $i=i_v$, $v\subseteq u^*$ then $f^*_i=(f^*\cup h_v)\comp H_{i,0}$,
\item if $i\notin\{i_v: v\subseteq u^*\}$ then $f^*_i\supseteq f^*\comp
H_{i,0}$ is such that $f^*_i(\alpha)=0$ for all $\alpha\in H_{i,0}[u^*]
\setminus u$. 
\end{itemize}
Suppose that $h:u\longrightarrow 2$ and let $f_i=f^*_i\cup (h\rest (u\cap
w^{p_i}))$. It should be clear that for each $i\leq N$ we have $f_i\in
F^{p_i}$ and $f_i\rest w^*=f_0\rest w^*$ (remember the choice of
$f^*$). Assume that $\{0,1\}\models \tau(f_0(\bar{\alpha}(0,0)),\ldots,f_0(
\bar{\alpha}(0,\ell)))=1$. Look at $v=h^{-1}[\{1\}]\cap u^*$ and the
corresponding $i_v$. By the above assumption and the choice of $h_v,f_{i_v}^*$
we have
\[\{0,1\}\models \tau(f_{i_v}(\bar{\alpha}(i_v,0)),\ldots,f_{i_v}(\bar{\alpha}
(i_v,\ell)))=1.\]
This shows that $\bigcup\limits_{i\leq N}f_i\in F^q$ and hence we conclude
$q\in\qlm$. If $|u^+\setminus u^*|\leq n/2$ then we proceed similarly:\quad
for $v\subseteq u^+\setminus u^*$ we choose distinct $i_v\in [1,N]$ such that
$w^{p_{i_v}}\cap u=u^*$. We pick $f^*$ as in the previous case and we define
$f^*_i:w^{p_i}\setminus u\longrightarrow 2$ (for $i\leq N$) as follows
\begin{itemize}
\item if $i=i_v$, $v\subseteq u^+\setminus u^*$ then $f^*_i=f^*\comp H_{i,0}$
and $(\forall\alpha\in u^+\setminus u^*)(f^*_i(H_{0,i}(\alpha))=1\
\Leftrightarrow\ \alpha\in v)$,
\item if $i\notin\{i_v: v\subseteq u^+\setminus u^*\}$ then $f^*_i\supseteq
f^*\comp H_{i,0}$ is such that $f^*_i(\alpha)=0$ for all $\alpha\in H_{i,0}
[u^+]\setminus u$. 
\end{itemize}
Next we argue like before to show that $q\in\qlm$. 

Checking that $q$ is a common upper bound of $p_0,\ldots,p_N$ is
straightforward. Finally, by the definition of $F^q$ and by \ref{1.2}(2) we
see that 
\[q\forces_{\qlm}\mbox{`` }\dot{a}_0\wedge \bigwedge_{j=1}^N (-\dot{a}_j)=0
\mbox{ ''.}\]
Thus we have proved that $\forces_{\qlm}$ ``$\tig^+_{1,N}(\dot{\bB})\leq
\mu^+$''. The same arguments show that $\forces_{\qlm}$ ``$\tig^+_{N,1}(
\dot{\bB})\leq\mu^+$'' (just considering $-\dot{a}_\alpha$ instead of
$\dot{a}_\alpha$ and $\{0,\ldots,N-1\}$, $\{N\}$ as the two groups of indexes
there). 

To show that the equalities hold one can prove even more: in $\V^{\qlm}$,
there is an independent subset of $\dot{\bB}$ of size $\mu$. The construction
of the set is easy once you note that if $p\in\qlm$,
$\alpha\in\lambda\setminus w^p$ and $w^q=w^p\cup\{\alpha\}$, $F^q=\{f\in
2^{\textstyle w^q}: f\rest w^p\in F^p\}$ then $q=(w^q,F^q)$ is a condition in
$\qlm$ stronger than $p$. \QED

\begin{conclusion}
\label{1.5}
Assume that $\mu^{<\mu}=\mu<\lambda\leq\chi$. Then there is a forcing notion
$\bP$ which does not change cardinalities and cofinalities and such that in
$\V^{\bP}$:\quad $2^\mu\geq\chi$ and there are Boolean algebras $\bB_0,\bB_1,
\bB_2,\ldots$ of size $\lambda$ satisfying
\[\ind^+_{n+1}(\bB_n)=\lambda^+\quad\mbox{ and }\quad\hd^+(\bB_n)=\hL^+(\bB_n)
=\ind^+(\bB_n)=\mu^+.\]
Consequently, in $\V^{\bP}$, for every non-principal ultrafilter $\cD$ on
$\omega$ we have 
\[\inv(\prod_{n<\omega}\bB_n/\cD)=\lambda^\omega\quad\mbox{ and }\quad
\prod_{n\in\omega}\inv(\bB_n)/\cD=\mu^\omega,\]
where $\inv\in\{\ind,\tig,\hd,\hL,s\}$.
\end{conclusion}

\Proof Let $\bP_0$ be the forcing notion adding $\chi$ many Cohen subsets of
$\mu$ (with conditions of size $<\mu$) and for $n>0$ let $\bP_n$ be $\qlm$. Let
$\bP$ be the ${<}\mu$--support product of the $\bP_n$'s (so if $\mu=\omega$
then $\bP$ is the finite support product of the $\bP_n$'s and otherwise it is
the full product). 

\begin{claim}
\label{1.5.1}
$\bP$ is a $\mu$--closed $\mu^+$--cc forcing notion of size $\chi^{<\mu}$.
\end{claim}

\noindent{\em Proof of the claim:}\qquad Modify the proof of \ref{1.3}.
\medskip

Let $\dot{\bB}_n$ be the $\bP_{n+1}$--name (and so $\bP$--name) for the
Boolean algebra added by forcing with $\bP_n$.

\begin{claim}
\label{1.5.2}
For $n\in\omega$, $\inv\in\{\ind,\tig,\hd,\hL,s\}$ we have
\[\forces_{\bP}\mbox{`` }\ind^+_{n+1}(\dot{\bB}_n)=\lambda^+\ \mbox{ and }\
\inv^+(\dot{\bB}_n)=\mu^+\mbox{ ''.}\]
\end{claim}

\noindent{\em Proof of the claim:}\qquad Repeat the proof of \ref{1.4} with
suitable changes to show that in $\V^{\bP}$, for each $n$, we have 
\[\ind^+_{n+1}(\dot{\bB}_n)=\lambda^+\quad\mbox{ and }\quad\tig^+_{1,2^n+n}(
\dot{\bB}_n)=\tig_{2^n+n,1}(\dot{\bB}_n)=\ind^+(\dot{\bB}_n)=\mu^+.\]
Now note that for a Boolean algebra $\bB$
\[\tig^{(+)}_{1,N}=\hd^{(+)}_N(\bB)\quad\mbox{ and }\quad \tig^{(+)}_{N,1}(
\bB)=\hL^{(+)}_N(\bB)\]
(and remember that $\ind^{(+)}(\bB)\leq s^{(+)}(\bB)\leq\hd^{(+)}(\bB),
\hL^{(+)}(\bB)$). 
\medskip

The ``consequently'' part of the conclusion should be clear (or see
\cite[Section 1]{RoSh:534}). \QED 

\begin{remark}
\label{1.6}
Note that the examples when the spread of ultraproduct is larger than the
ultraproduct of the spreads which were known before provided ``a successor''
difference only. Conclusion \ref{1.5} shows that the jump can be larger, but
we do not know if one can get it in ZFC (i.e.~assuming suitable cardinal
arithmetic only).
\end{remark}

\begin{problem}
\label{1.7}
Can one improve \ref{1.4} getting it for $N=n+1$?
\end{problem}

\section{Irredundance of products}
In theorem \ref{2.1} below we answer \cite[Problem 24]{M2}. A parallel
question for free products of Boolean algebras will be addressed in the
next section. It should be noted here that the proof of the ZFC result was
written as a result of an analysis why a forcing proof of consistency of an
inequality (similar to the one from the next section) failed.

\begin{theorem}
\label{2.1}
For Boolean algebras $\bB_0,\bB_1$:
\[\irr(\bB_0\times\bB_1)=\max\{\irr(\bB_0),\irr(\bB_1)\}.\]
\end{theorem}

\Proof Clearly $\irr(\bB_0\times\bB_1)\geq\max\{\irr(\bB_0),\irr(\bB_1)\}$, so
we have to deal with the converse inequality only. Assume that a sequence
$\bar{x}=\langle (x^0_\alpha,x^1_\alpha):\alpha<\lambda\rangle\subseteq
\bB_0\times\bB_1$ is irredundant. Thus, for each $\alpha<\lambda$, we have
homomorphisms $f^0_\alpha,f^1_\alpha:\bB_0\times\bB_1\longrightarrow\{0,1\}$
such that $f^0_\alpha(x^0_\alpha,x^1_\alpha)=0$, $f^1_\alpha(x^0_\alpha,
x^1_\alpha)=1$ and
\[(\forall\beta\in\lambda\setminus\{\alpha\})(f^0_\alpha(x^0_\beta,
x^1_\beta)=f^1_\alpha(x^0_\beta,x^1_\beta)).\]
By shrinking the sequence $\bar{x}$ if necessary, we may assume that one
of the following occurs:
\begin{enumerate}
\item[(i)]\ \ \quad $(\forall\alpha<\lambda)(f^0_\alpha(1,0)=f^1_\alpha(1,0)
=0)$,
\item[(ii)]\ \quad $(\forall\alpha<\lambda)(f^0_\alpha(1,0)=f^1_\alpha(1,0)
=1)$,
\item[(iii)] \quad $(\forall\alpha<\lambda)(f^0_\alpha(1,0)=0\ \&\
f^1_\alpha(1,0)=1)$,
\item[(iv)]\ \quad $(\forall\alpha<\lambda)(f^0_\alpha(1,0)=1\ \&\
f^1_\alpha(1,0)=0)$.
\end{enumerate}
If the first clause occurs then we may define (for $\alpha<\lambda$)
homomorphisms $h^0_\alpha,h^1_\alpha:\bB_1\longrightarrow\{0,1\}$ by
$h^\ell_\alpha(x)=f^\ell_\alpha(1,x)$ (remember that in this case we
have $f^\ell_\alpha(0,1)=1$). Clearly these homomorphisms witness that
the sequence $\langle x^1_\alpha:\alpha<\lambda\rangle\subseteq\bB_1$ is
irredundant (and thus $\irr^+(\bB_1)>\lambda$). Similarly, if (ii) holds
then the sequence $\langle x^0_\alpha:\alpha<\lambda\rangle\subseteq
\bB_0$ is irredundant and $\irr^+(\bB_0)>\lambda$.

Since $f^\ell_\alpha(1,0)=0\ \Leftrightarrow\ f^\ell_\alpha(0,1)=1$ and
the algebras $\bB_0,\bB_1$ are in symmetric positions, we may assume
that clause (iv) holds, so $f^\ell_\alpha(0,1)=\ell$ (for $\ell<2$,
$\alpha<\lambda$).

For $\alpha<\lambda$ and $\ell<2$ let $g^\ell_\alpha:\lambda\longrightarrow
2$ be given by $g^\ell_\alpha(\beta)=f^\ell_\alpha(x^0_\beta,x^1_\beta)$
for $\beta<\lambda$. Note that $\beta\neq\alpha$ implies $g^0_\alpha(\beta)=
g^1_\alpha(\beta)$ (remember the choice of the $f^\ell_\alpha$'s). Next, for
$\ell<2$ let $F_\ell=\{g^\ell_\alpha:\alpha<\lambda\}$ and let $\bB^*_\ell$
be the algebra $\bB_{(\lambda,F_\ell)}$ (see \ref{1.1}(1)).

\begin{claim}
\label{2.1.1}
Assume that $A\subseteq\lambda$ and $\ell<2$ are such that
\begin{enumerate}
\item[$(\boxtimes^\ell_A)$] the mappings $\{x_\beta:\beta\in A\}
\longrightarrow \{0,1\}: x_\beta\mapsto g^k_\alpha(\beta)$ (for $k=0,1$ and
$\alpha\in A$) extend to homomorphisms from $\langle x_\beta:\beta\in A
\rangle_{\bB^*_\ell}$ onto $\{0,1\}$.
\end{enumerate}
Then the sequence $\langle x_\alpha^\ell:\alpha\in A\rangle\subseteq\bB_\ell$
is irredundant.
\end{claim}

\noindent{\em Proof of the claim:}\qquad First note that the assumption
$(\boxtimes^\ell_A)$ implies that the sequence $\langle x_\beta:\beta\in
A\rangle\subseteq\bB^*_\ell$ is irredundant. Now, the mapping $x_\beta^\ell
\mapsto x_\beta$ extends to a homomorphism from the algebra $\langle
x^\ell_\beta: \beta<\lambda\rangle_{\bB_\ell}$ onto $\bB^*_\ell$. [Why? Note
that, since $f^0_\alpha(1,0)=1=f^1_\alpha(0,1)$, the mappings $x^\ell_\beta
\mapsto f^\ell_\alpha(x^0_\beta,x^1_\beta)=g^\ell_\alpha(\beta)$ extend to
homomorphisms from $\bB_\ell$ onto $\{0,1\}$. Now look at the definition of
the algebra $\bB^*_\ell$; remember \ref{1.2}(2).] Consequently we get that the
sequence $\langle x^\ell_\beta:\beta\in A\rangle\subseteq\bB_\ell$ is
irredundant. 
\medskip

It follows from claim \ref{2.1.1} that if there are $A\in [\lambda]^{
\textstyle\lambda}$ and $\ell<2$ such that $(\boxtimes^\ell_A)$ holds
true then the algebra $\bB_\ell$ has an irredundant sequence of length
$\lambda$ (i.e.~$\irr^+(\bB_\ell)>\lambda$). So the proof of the theorem
will be concluded when we show the following claim.

\begin{claim}
\label{2.1.2}
Let $\ell<2$. Assume that there is no $A\in [\lambda]^{\textstyle\lambda}$
such that $(\boxtimes^\ell_A)$ holds. Then $s^+(\bB_{1-\ell})>\lambda$ (so
$\irr^+(\bB_{1-\ell})>\lambda$ too).
\end{claim}

\noindent{\em Proof of the claim:}\qquad By induction on $\xi<\lambda$
we build a sequence $\langle (u_\xi,v_\xi):\xi<\lambda\rangle$ such that for
each $\xi<\lambda$:
\begin{enumerate}
\item[(a)] $u_\xi,v_\xi\in [\lambda]^{\textstyle <\omega}$ are disjoint,
\item[(b)] $(u_\xi\cup v_\xi)\cap \bigcup\limits_{\zeta<\xi} (u_\zeta\cup
v_\zeta)=\emptyset$,
\item[(c)] $\bB^*_\ell\models\bigwedge\limits_{\gamma\in u_\xi}x_\gamma\wedge
\bigwedge\limits_{\gamma\in v_\xi}(-x_\gamma)= 0$,
\item[(d)] $\bB^*_{1-\ell}\models\bigwedge\limits_{\gamma\in u_\xi}x_\gamma
\wedge\bigwedge\limits_{\gamma\in v_\xi}(-x_\gamma)\neq 0$.
\end{enumerate}
Suppose we have defined $u_\zeta,v_\zeta$ for $\zeta<\xi$. The set
$A=\lambda\setminus\bigcup\limits_{\zeta<\xi}(u_\zeta\cup v_\zeta)$ is
of size $\lambda$, so (by our assumptions) $(\boxtimes^\ell_A)$ fails.
This means that one of the mappings
\[\{x_\beta:\beta\in A\}\longrightarrow\{0,1\}:x_\beta\mapsto g^k_\alpha(
\beta),\qquad(k=0,1,\ \alpha\in A)\]
does not extend to a homomorphism from $\langle x_\beta:\beta\in A\rangle_{
\bB^*_\ell}$. But, by the definition of $\bB^*_\ell$, the mappings $x_\beta
\mapsto g^\ell_\alpha(\beta)$ do extend (see \ref{1.2}(1)). So we find finite
disjoint sets $u_\xi,v_\xi\subseteq A$ such that $\bB^*_\ell\models\bigwedge
\limits_{\gamma\in u_\xi}x_\gamma\wedge\bigwedge\limits_{\gamma\in v_\xi}(-
x_\gamma)= 0$, but for some $\alpha<\lambda$, $g^{1-\ell}_\alpha\rest u_\xi
\equiv 1$ and $g^{1-\ell}_\alpha\rest v_\xi\equiv 0$. The latter implies that
$\bB^*_{1-\ell}\models\bigwedge\limits_{\gamma\in u_\xi}x_\gamma\wedge\bigwedge
\limits_{\gamma\in v_\xi}(-x_\gamma)\neq 0$. This finishes the construction.

The demand (d) means that (by \ref{1.2}) for each $\xi<\lambda$ we find
$\alpha_\xi<\lambda$ such that $g^{1-\ell}_{\alpha_\xi}\rest u_\xi\equiv 1$
and $g^{1-\ell}_{\alpha_\xi}\rest v_\xi\equiv 0$. On the other hand, by (c),
there is no $\alpha<\lambda$ such that $g^\ell_\alpha\rest u_\xi\equiv 1$ and
$g^\ell_\alpha\rest v_\xi\equiv 0$. But now, if $\alpha\notin u_\xi\cup v_\xi$
then $g^{1-\ell}_\alpha\rest (u_\xi\cup v_\xi)=g^\ell_\alpha\rest (u_\xi\cup
v_\xi)$, so necessarily $\alpha_\xi\in u_\xi\cup v_\xi$. Let $y_\xi=\bigwedge
\limits_{\gamma\in u_\xi}x^{1-\ell}_\gamma\wedge\bigwedge\limits_{\gamma\in
v_\xi}(-x^{1-\ell}_\gamma)\in\bB_{1-\ell}$ and $h_\xi:\langle x^{1-
\ell}_\beta:\beta<\lambda\rangle_{\bB_{1-\ell}}\longrightarrow \{0,1\}$ be a
homomorphism defined by $h_\xi(x^{1-\ell}_\beta)=f^{1-\ell}_{\alpha_\xi}(
x^0_\beta,x^1_\beta)=g^{1-\ell}_{\alpha_\xi}(\beta)$. It follows from
the above discussion that ($h_\xi$ is well defined and)
\[h_\xi(y_\zeta)=1\quad\mbox{ if and only if }\quad \xi=\zeta,\]
showing that the sequence $\langle y_\xi:\xi<\lambda\rangle$ is ideal 
independent (and irredundant). This finishes the proof of the claim and that 
of the theorem. \QED

\section{Forcing for spread and irredundance}
In this section we show that, consistently, there is a Boolean algebra
$\bB$ such that $\irr(\bB)<s(\bB\circledast\bB)$. This gives a partial answer
to \cite[Problem 27]{M2}. Moreover, it shows that a statement parallel to 
\ref{2.1} for the free product (instead of product) is not provable in ZFC.
Note that before trying to answer \cite[Problem 27]{M2} in ZFC one should
first construct a ZFC example of a Boolean algebra $\bB$ such that $\irr(\bB)
<|\bB|$ -- so far no such example is known.

\begin{definition}
\label{3.1}
\begin{enumerate}
\item We define a forcing notion $\bQ^*$ by:

\noindent{\bf a condition} is a tuple $p=\langle u^p,\langle f^p_{0,\alpha},
f^p_{1,\alpha}, f^p_{2,\alpha}:\alpha\in u^p\rangle\rangle$ such that
\begin{enumerate}
\item[(a)] $u^p\subseteq\omega_1$ is finite,
\item[(b)] $f^p_{\ell,\alpha}:u^p\times 2\longrightarrow\{0,1\}$ for
$\ell<3,\alpha\in u^p$,
\item[(c)] $f^p_{0,\alpha}\rest (u^p\cap\alpha)\times 2=f^p_{1,\alpha}\rest
(u^p\cap\alpha)\times 2=f^p_{2,\alpha}\rest (u^p\cap\alpha)\times 2$ for
$\alpha\in u^p$,
\item[(d)] $f^p_{0,\alpha}(\alpha,0)=1$, $f^p_{0,\alpha}(\alpha,1)=0$
(for $\alpha\in u^p$),
\item[(e)] $f^p_{1,\alpha}(\alpha,0)=0$, $f^p_{1,\alpha}(\alpha,1)=1$
(for $\alpha\in u^p$),
\item[(f)] $f^p_{0,\alpha}(\beta,0)=0$ or  $f^p_{1,\alpha}(\beta,1)=0$
(for distinct $\alpha,\beta\in u^p$),
\item[(g)] $f^p_{0,\alpha}(\beta,0)=0$ or $f^p_{2,\alpha}(\beta,1)=0$
(for $\alpha,\beta\in u^p$),
\item[(h)] $f^p_{1,\alpha}(\beta,1)=0$ or $f^p_{2,\alpha}(\beta,0)=0$
(for $\alpha,\beta\in u^p$),
\item[(i)] $f^p_{2,\alpha}(\beta,0)=0$ or $f^p_{2,\alpha}(\beta,1)=0$
(for $\alpha,\beta\in u^p$);
\end{enumerate}
{\bf the order} is defined by:\quad $p\leq q$\quad if and only if\quad
$u^p\subseteq u^q$, and $f^q_{\ell,\alpha}\rest (u^p\times 2)=f^p_{\ell,  
\alpha}$ for $\alpha\in u^p$, $\ell<3$ and for each $\alpha\in u^q$,
$\ell<3$:
\[f^q_{\ell,\alpha}\rest (u^p\times 2)\in\{f^p_{k,\beta}:\beta\in u^p,\
k<3\}.\]
\item For a condition $p\in\bQ^*$ let $\bB^*_p$ be the algebra $\bB_{(w,F)}$,
where $w=u^p\times 2$ and $F=\{f^p_{\ell,\alpha}:\alpha\in u^p,\ell<3\}$ (see
\ref{1.1}(1)).
\item Let $\dot{\bB}^*$, $\dot{f}_{\ell,\alpha}$ (for $\ell<3$, $\alpha<
\omega_1$) be $\bQ^*$-names such that
\[\forces_{\bQ^*}\mbox{`` }\dot{\bB}^*=\bigcup\{\bB_p^*:p\in\Gamma_{\bQ^*}\}\
\mbox{ and }\ \dot{f}_{\ell,\alpha}=\bigcup\{f^p_{\ell,\alpha}: p\in \Gamma_{
\bQ^*}, \alpha\in u^p\}\mbox{ ''.}\]
\end{enumerate}
\end{definition}

\begin{proposition}
\label{3.2}
\begin{enumerate}
\item $\bQ^*$ is a ccc forcing notion.
\item If $p,q\in\bQ^*$, $p\leq q$ then $\bB_p$ is a subalgebra of $\bB_q$.
\item In $\V^{\bQ^*}$, $\dot{f}_{\ell,\alpha}:\omega_1\times 2\longrightarrow
2$ (for $\alpha<\omega_1$ and $\ell<3$) and $\dot{\bB}^*$ is the Boolean
algebra $\bB_{(w,F)}$, where $w=\omega_1\times 2$ and $F=\{\dot{f}_{\ell,
\alpha}:\alpha<\omega_1,\ell<3\}$.
\end{enumerate}
\end{proposition}

\Proof 1)\ \ \ Suppose that $\cA\subseteq\bQ^*$ is uncountable. Applying
$\Delta$--system arguments find $p,q\in\cA$ such that letting $u^*=u^p\cap
u^q$ we have:
\begin{enumerate}
\item[(i)]\ \ \ $\max(u^*)<\min(u^p\setminus u^*)\leq\max(u^p\setminus u^*)<
\min(u^q\setminus u^*)$,
\item[(ii)]\ \ $|u^p|=|u^q|$ and if $H:u^p\longrightarrow u^q$ is the order
isomorphism, $\alpha\in u^p$ and $\ell<3$ then $f^p_{\ell,\alpha}=f^q_{\ell,H(
\alpha)}\comp(H\times{\rm id})$.
\end{enumerate}
Now let $u^r=u^p\cup u^q$ and for $\ell<3$ and $\alpha\in u^r$ let:
\[f^r_{\ell,\alpha}=\left\{\begin{array}{ll}
f^p_{\ell,\alpha}\cup f^q_{\ell,\alpha} &\mbox{ if }\alpha\in u^p\cap u^q,\\
f^p_{\ell,\alpha}\cup f^q_{2,H(\alpha)}\rest(u^q\setminus u^p) &\mbox{ if }
\alpha\in u^p\setminus u^q,\\
f^q_{\ell,\alpha}\cup f^p_{2,H^{-1}(\alpha)}\rest(u^p\setminus u^q)&\mbox{ if }
\alpha\in u^q\setminus u^p.
			   \end{array}\right.\]
It is a routine to check that this defines a condition in $\bQ^*$ stronger
than both $p$ and $q$.
\medskip

\noindent 2)\ \ \ Should be clear.
\medskip

\noindent 3)\ \ \ Note that if $p\in\bQ^*$, $\alpha_0\in u^p$ and $\beta\in
\omega_1\setminus u^p$ then letting $u^q=u^p\cup\{\beta\}$ and
\[f^q_{\ell,\alpha}=\left\{\begin{array}{ll}
f^p_{\ell,\alpha}\cup\{((\beta,0),0),((\beta,1),0)\}&\mbox{ if }\alpha\in
u^p,\ \ell<3,\\
f^p_{2,\alpha_0}\cup\{((\beta,0),1-\ell),((\beta,1),\ell)\}&\mbox{ if }\alpha
=\beta,\ \ell<2,\\
f^p_{2,\alpha_0}\cup\{((\beta,0),0),((\beta,1),0)\}&\mbox{ if }\alpha=\beta,\
\ell=2,
\end{array}\right.\]
we get a condition $q\in\bQ^*$ stronger than $p$ and such that $\beta\in
w^q$. Now, the rest should be clear. \QED

\begin{proposition}
\label{3.3}
$\forces_{\bQ^*}$`` $s(\dot{\bB}^*\circledast\dot{\bB}^*)=\omega_1$ ''.
\end{proposition}

\Proof To avoid confusion between the two copies of $\dot{\bB}^*$ in
$\dot{\bB}^*\circledast\dot{\bB}^*$, let us denote an element $a\wedge
b\in\dot{\bB}^*\circledast\dot{\bB}^*$ such that $a$ is from the first
copy of $\dot{\bB}^*$ and $b$ is from the second one, by $\langle
a,b\rangle$. With this convention, for each $\alpha<\omega_1$ let
$\dot{y}_\alpha=\langle x_{\alpha,0},x_{\alpha,1}\rangle$ and let
$\dot{f}_\alpha:\dot{\bB}^*\circledast\dot{\bB}^*\longrightarrow\{0,1\}$
be a homomorphism such that (for $\beta<\omega_1$, $i<2$)
\[\dot{f}_\alpha(\langle x_{\beta,i},1\rangle)=\dot{f}_{0,\alpha}(\beta,i)
\quad\mbox{and}\quad\dot{f}_\alpha(\langle 1,x_{\beta,i}\rangle)=\dot{f}_{1,
\alpha}(\beta,i).\]
Note that, by \ref{3.1}(d,e), for each $\alpha<\omega_1$
\[\dot{f}_\alpha(\dot{y}_\alpha)=\dot{f}_{0,\alpha}(\alpha,0)\wedge
\dot{f}_{1,\alpha}(\alpha,1)=1,\]
and if $\beta\in\omega_1\setminus\{\alpha\}$ then (by \ref{3.1}(f))
\[\dot{f}_\alpha(\dot{y}_\beta)=\dot{f}_{0,\alpha}(\beta,0)\wedge
\dot{f}_{1,\alpha}(\beta,1)=0.\]
Hence we conclude that
\[\forces_{\bQ^*}\mbox{`` }\langle\dot{y}_\alpha:\alpha<\omega_1\rangle
\mbox{ is ideal--independent '',}\]
finishing the proof. \QED

\begin{theorem}
\label{3.4}
$\forces_{\bQ^*}$ `` $\irr^+_5(\dot{\bB}^*)=\omega_0$ ''.
\end{theorem}

\Proof Let $\langle\dot{a}_\beta:\beta<\omega_1\rangle$ be a $\bQ^*$--name for
an $\omega_1$--sequence of elements of $\dot{\bB}^*$, $p\in\bQ^*$. For $\beta<
\omega_1$ choose a condition $p_\beta\geq p$, a Boolean term $\tau_\beta$,
ordinals $\bar{\alpha}(\beta,0)\leq\ldots\leq\bar{\alpha}(\beta,\ell_\beta)<
\omega_1$ and $\bar{\imath}(\beta,0),\ldots,\bar{\imath}(\beta,\ell_\beta)\in
\{0,1\}$ such that
\[p_\beta\forces_{\bQ^*}\dot{a}_\beta=\tau_\beta(x_{\bar{\alpha}(\beta,0),
\bar{\imath}(\beta,0)},\ldots,x_{\bar{\alpha}(\beta,\ell_\beta),\bar{\imath}
(\beta,\ell_\beta)}).\]
Applying standard ``cleaning procedure'' we may assume that for
$\beta,\beta_0,\beta_1<\omega_1$:
\begin{enumerate}
\item[(i)]\ \ \ $\tau_\beta=\tau$, $\ell_\beta=\ell$,
\item[(ii)]\ \  $\{(\bar{\alpha}(\beta,j),\bar{\imath}(\beta,j)):
j\leq\ell\}=u^{p_\beta}\times 2$ is an enumeration which does not depend
on $\beta$ if we treat it modulo $\otp$ (so $2\cdot |u^{p_\beta}|=\ell+1$ and
we may write $\tau(x_{\gamma,i}:\gamma\in u^{p_\beta}, i<2)$),
\item[(iii)]\   $\{u^{p_\beta}:\beta<\omega_1\}$ forms a $\Delta$--system
of sets with the heart $u^*$, and if $\beta_0<\beta_1<\omega_1$ then
\[\max(u^*)<\min(u^{p_{\beta_0}}\setminus u^*)\leq\max(u^{p_{\beta_0}}
\setminus u^*)<\min(u^{p_{\beta_1}}\setminus u^*),\]
\item[(iv)]\ \  $|u^{p_{\beta_0}}|=|u^{p_{\beta_1}}|$ and if $H_{\beta_0,
\beta_1}:u^{p_{\beta_0}}\longrightarrow u^{p_{\beta_1}}$ is the order
preserving mapping then $f^{p_{\beta_0}}_{k,\alpha}=f^{p_{\beta_1}}_{k,H(
\alpha)}\comp (H_{\beta_0,\beta_1}\times{\rm id})$ (for $\alpha\in u^{p_{
\beta_0}}$, $k<3$).
\end{enumerate}
Now we are going to define a condition $q$ stronger than $p_0,\ldots,p_5$.
We put $u^q=\bigcup\limits_{i<6}u^{p_i}$ and we define functions
$f^q_{\ell,\alpha}:u^q\times 2\longrightarrow 2$ (for $\alpha\in u^q$
and $\ell<3$) as follows.
\begin{enumerate}
\item[$(\boxtimes)$]  If $\alpha\in u^*$, $\ell<3$ then $f^q_{\ell,\alpha}=
\bigcup\limits_{i<6}f^{p_i}_{\ell,\alpha}$.
\item[$(\boxplus_0)$] If $\alpha\in u^{p_0}\setminus u^*$ then
\[\begin{array}{ll}
f^q_{0,\alpha}=& f^{p_0}_{0,\alpha}\cup f^{p_1}_{0,H_{0,1}(\alpha)}\cup
f^{p_2}_{0,H_{0,2}(\alpha)}\cup f^{p_3}_{2,H_{0,3}(\alpha)}\cup
f^{p_4}_{2,H_{0,4}(\alpha)}\cup f^{p_5}_{2,H_{0,5}(\alpha)},\\
f^q_{1,\alpha}=& f^{p_0}_{1,\alpha}\cup f^{p_1}_{2,H_{0,1}(\alpha)}\cup
f^{p_2}_{2,H_{0,2}(\alpha)}\cup f^{p_3}_{1,H_{0,3}(\alpha)}\cup
f^{p_4}_{1,H_{0,4}(\alpha)}\cup f^{p_5}_{2,H_{0,5}(\alpha)},\\
f^q_{2,\alpha}=& \bigcup\limits_{i<6} f^{p_i}_{2,H_{0,i}(\alpha)}.
\end{array}\]
\item[$(\boxplus_i)$] If $\alpha\in u^{p_i}\setminus u^*$, $0<i<6$ and
$\ell<3$ then $f^q_{\ell,\alpha}=f^{p_i}_{\ell,\alpha}\cup
\bigcup\limits_{j\neq i} f^{p_j}_{2,H_{i,j}(\alpha)}$.
\end{enumerate}
It follows from {\bf (iv)} and \ref{3.1}(c) that the functions $f^q_{\ell,
\alpha}$ are well defined.

\begin{claim}
\label{3.3.1}
The tuple $q=\langle u^q,\langle f^q_{\ell,\alpha}: \ell<3,\alpha\in
u^q\rangle\rangle$ is a condition in $\bQ^*$ stronger than $p_0,\ldots,p_5$.
\end{claim}

\noindent{\em Proof of the claim:}\qquad To show that $q\in\bQ^*$ one
has to check the demands (a)--(i) of \ref{3.1}. The only possible
problems could be caused by clauses (f)--(i). If functions
$f^q_{\ell,\alpha}$  were defined in clauses $(\boxtimes)$, $(\boxplus_i)$
then easily these demands are met. To deal with instances of
$(\boxplus_0)$ (i.e.~when $\alpha\in u^{p_0}\setminus u^*$) note that
in the definition of $f^q_{\ell,\alpha}$ ($\ell<2$, $\alpha\in u^{p_0}
\setminus u^*$) a part of the form $f^{p_j}_{\ell, H_{i,j}(\alpha)}$
``meets'' $f^{p_j}_{2,H_{i,j}(\alpha)}$ on the side of
$f^q_{1-\ell,\alpha}$. Therefore, by (g), (h) of \ref{3.1}, we have no
problems with checking demand (f). Clause \ref{3.1}(i) is immediate and
(g), (h) should be clear too.

\begin{claim}
\label{3.3.2}
\[\tau(x_{\gamma,i}:\gamma\in u^{p_0}, i<2)\in\langle\tau(x_{\gamma,i}:
\gamma\in u^{p_j}, i<2): 0<j<6\rangle_{\bB^*_q}.\]
Consequently\quad $q\forces_{\bQ^*}$`` $\dot{a}_0\in\langle \dot{a}_j: 0<j<6
\rangle_{\dot{\bB}^*}$''.
\end{claim}

\noindent{\em Proof of the claim:}\qquad Suppose that
\[\tau(x_{\gamma,i}: \gamma\in u^{p_0}, i<2)\notin\langle\tau(x_{\gamma,i}:
\gamma\in u^{p_j}, i<2): 0<j<6\rangle_{\bB^*_q}.\]
Then we find two homomorphisms $h_0,h_1:\bB^*_q\longrightarrow\{0,1\}$
such that
\[\begin{array}{l}
h_0(\tau(x_{\gamma,i}:\gamma\in u^{p_0}, i<2))\neq h_1(\tau(x_{\gamma,i}:
\gamma\in u^{p_0},i<2))\ \mbox{ but}\\
h_0(\tau(x_{\gamma,i}:\gamma\in u^{p_j},i<2))=h_1(\tau(x_{\gamma,i}:\gamma\in
u^{p_j},i<2))\ \mbox{ for }0<j<6.\\
\end{array}\]
By the definition of the algebra $\bB^*_q$ each its homomorphism into $\{0,1\}$
is generated by one of the functions $f^q_{\ell,\alpha}$ (for $\ell<3$, $\alpha
\in u^q$). So we find $\ell_0,\ell_1<3$ and $\alpha_0,\alpha_1\in u^q$ such
that $h_k\supseteq f^q_{\ell_k,\alpha_k}$. Now we have to consider several
cases corresponding to the way the $f^q_{\ell_k,\alpha_k}$ were defined.
\medskip

\noindent{\sc Case A:}\qquad $\alpha_k\in u^*$, $\alpha_{1-k}\in
u^{p_i}$, $i<6$.\\
Then look at the definition $(\boxtimes)$ of $f^q_{\ell_k,\alpha_k}$ --
it copies $f^{p_0}_{\ell_k,\alpha_k}$ everywhere (remember {\bf (iv)}).
On the other hand, whatever clause was used to define $f^q_{\ell_{1-k},
\alpha_{1-k}}$, there is $j\in (0,6)$ such that $f^q_{\ell_{1-k},\alpha_{1-
k}}\rest (u^{p_j}\times 2)$ is a copy of $f^q_{\ell_{1-k},\alpha_{1-k}}\rest
(u^{p_0}\times 2)$. Hence we may conclude that (for this $j$)
\[h_{1-k}(\tau(x_{\gamma,i}:\gamma\in u^{p_j},i<2))\neq h_k(\tau(x_{\gamma,i}:
\gamma\in u^{p_j},i<2)),\]
a contradiction.
\medskip

\noindent{\sc Case B:}\qquad $\alpha_k\in u^{p_0}\setminus u^*$, $\alpha_{1-k}
\in u^{p_i}\setminus u^*$, $0<i<6$.\\
Then we repeat the argument of the previous Case, choosing $j$ in such a
way that $j\neq i$ and:\quad if $\ell_k=0$ then $j\in\{1,2\}$,\quad if $\ell_k
=1$ then $j\in\{3,4\}$.
\medskip

\noindent{\sc Case C:}\qquad $\alpha_k\in u^{p_{i'}}\setminus u^*$,
$\alpha_{1-k}\in u^{p_{i''}}\setminus u^*$, $0<i',i''<6$.\\
Like above, but now take $j\in\{1,\ldots,5\}\setminus\{i',i''\}$.
\medskip

\noindent{\sc Case D:}\qquad $\alpha_0,\alpha_1\in u^{p_0}\setminus
u^*$.\\
This is the most complicated case. We may repeat the previous argument
in some cases letting:
\[j=\left\{\begin{array}{ll}
1 &\mbox{ if }\ (\ell_0,\ell_1)\in\{(0,0),(0,2),(2,0),(2,2)\}\\
3 &\mbox{ if }\ (\ell_0,\ell_1)\in\{(1,1),(1,2),(2,1)\}\\
\end{array}\right.\]
This leaves us with two symmetrical cases: $(\ell_0,\ell_1)=(0,1)$ or
$(\ell_0,\ell_1)=(1,0)$. So suppose that $\ell_0=0$, $\ell_1=1$ and let
\[x\stackrel{\rm def}{=}h_0(\tau(x_{\gamma,i}:\gamma\in u^{p_5},i<2))=
h_1(\tau(x_{\gamma,i}:\gamma\in u^{p_5},i<2)).\]
Since $f^q_{0,\alpha_0}\rest (u^{p_4}\times 2)$ is a copy of $f^q_{0,\alpha_0}
\rest (u^{p_5}\times 2)$ we conclude that
\[x=h_0(\tau(x_{\gamma,i}:\gamma\in u^{p_4},i<2))=h_1(\tau(x_{\gamma,i}:\gamma
\in u^{p_4},i<2)),\]
and, since $f^q_{1,\alpha_1}\rest (u^{p_4}\times 2)$ is a copy of $f^q_{1,
\alpha_1}\rest (u^{p_0}\times 2)$ we get
\begin{enumerate}
\item[$(\square)$]\quad $x=h_1(\tau(x_{\gamma,i}:\gamma\in u^{p_0},i<2))$.
\end{enumerate}
Next, $f^q_{1,\alpha_1}\rest (u^{p_2}\times 2)$ is a copy of $f^q_{1,\alpha_1}
\rest (u^{p_5}\times 2)$ and therefore
\[x=h_1(\tau(x_{\gamma,i}:\gamma\in u^{p_2},i<2))=h_0(\tau(x_{\gamma,i}:\gamma
\in u^{p_2},i<2)).\]
But $f^q_{0,\alpha_0}\rest (u^{p_2}\times 2)$ is a copy of $f^q_{0,\alpha_0}
\rest (u^{p_0}\times 2)$, so we conclude that
\begin{enumerate}
\item[$(\circledcirc)$]\quad $x=h_0(\tau(x_{\gamma,i}:\gamma\in u^{p_0},i<2))$.
\end{enumerate}
But now $(\square)+(\circledcirc)$ contradict the choice of $h_0,h_1$.
The other case is similar.
\medskip

\noindent This finishes the proof of the claim and of the theorem. \QED

\begin{conclusion}
\label{3.5}
It is consistent that there exists a Boolean algebra $\bB$ such that
\[\omega_0=\irr(\bB)\quad\mbox{ and }\quad s(\bB\circledast\bB)=\irr(\bB
\circledast\bB)=\omega_1.\qquad\QED\]
\end{conclusion}

\begin{remark}
We may use any cardinal $\mu=\mu^{<\mu}$ instead of $\omega$ and
$\mu^+$ instead of $\omega_1$ in \ref{3.1} and then \ref{3.2},
\ref{3.3}. But we do not know if the difference between the respective 
cardinal invariants can be larger.
\end{remark}

\begin{problem}
Is it consistent that there is a Boolean algebra $\bB$ such that

\qquad $(\irr(\bB))^+<|\bB|$ ?\qquad $(\irr(\bB))^+<s(\bB\circledast\bB)$ ?
\end{problem}

\section{Forcing a superatomic Boolean algebra}
In this section we give partial answers to \cite[Problems 73, 77, 78]{M2}
showing that, consistently, there is a superatomic Boolean algebra $\bB$ such
that $s(\bB)=\inc(\bB)<\irr(\bB)=\Id(\bB)<\Sub(\bB)$. The forcing notion we
use is a variant of the one of Martinez \cite{Ma92}, which in turn was a
modification of the forcing notion used in Baumgartner Shelah \cite{BaSh:254}.
For more information on superatomic Boolean algebras we refer the reader to
Koppelberg \cite{Ko89}, Roitman \cite{Rt89} and Monk \cite{M2}. 

\begin{definition}
\label{4.1}
Let $\kappa$ be a cardinal. For a pair $s=(\alpha,\xi)\in\kappa^+\times\kappa$
we will write $\alpha(s)=\alpha$ and $\xi(s)=\xi$. We define a forcing notion
$\bPk$ as follows: 
\smallskip

\noindent {\bf a condition}\quad is a tuple 
\[p=\left\langle w^p,u^p,a^p,\langle f^p_s:s\in u^p\rangle,\langle y^p_{s_0,
s_1}:s_0,s_1\in u^p,\ s_0\neq s_1,\ \alpha(s_0)\leq\alpha(s_1)\rangle\right
\rangle\] 
such that 
\begin{enumerate}
\item[(a)] $a^p\subseteq w^p\in [\kappa^+]^{\textstyle <\!\kappa}$, $u^p\in
[w^p\times\kappa]^{\textstyle<\!\kappa}$, and $\alpha\in w^p\ \Rightarrow\
(\alpha,0),(\alpha,1)\in u^p$,  
\item[(b)] for $s\in u^p$, $f^p_s:u^p\longrightarrow\{0,1\}$ is such that
$f^p_s(s)=1$ and 
\[(\forall t\in u^p)(\alpha(t)\leq\alpha(s)\ \&\ t\neq s\ \ \Rightarrow\ \
f^p_s(t)=0),\]  
\item[(c)] if $\alpha<\beta$, $\alpha,\beta\in a^p$ then $f^p_{\alpha,0}(
\beta,0)=f^p_{\alpha,1}(\beta,0)$,
\item[(d)] if $s_0,s_1\in u^p$ are distinct, $\alpha(s_0)\leq\alpha(s_1)$ then

$y_{s_0,s_1}^p\in [u^p\cap (\alpha(s_0)\times\kappa)]^{\textstyle<\!\omega}$
and for every $t\in u^p$
\[f^p_t(s_0)=1\ \&\ f^p_t(s_1)\neq f^p_{s_0}(s_1)\quad\Rightarrow\quad
(\exists s\in y^p_{s_0,s_1})(f^p_t(s)=1);\]
\end{enumerate}

\noindent {\bf the order}\quad is given by\quad $p\leq q$ if and only if
$w^p\subseteq w^q$, $u^p\subseteq u^q$, $a^p=a^q\cap w^p$, $y^q_{s_0,s_1}=
y^p_{s_0,s_1}$ (for distinct $s_0,s_1\in u^p$ such that $\alpha(s_0)\leq\alpha(
s_1)$), $f^p_s\subseteq f^q_s$ (for $s\in u^p$) and 
\[(\forall s\in u^q)(\exists t\in u^p)(f^q_s\rest u^p=f^p_t\ \mbox{ or }\
f^q_s\rest u^p={\bf 0}_{u^p} ).\]
\end{definition}

\begin{definition}
\label{isomorph}
We say that conditions $p,q\in\bPk$ are {\em isomorphic} if there is a
bijection $H:u^p\longrightarrow u^q$ (called {\em the isomorphism from $p$ to
$q$}) such that 
\begin{enumerate}
\item $(\forall s\in u^p)(\otp(\alpha(s)\cap w^p)=\otp(\alpha(H(s))\cap w^q)\ \
\&\ \ \xi(s)=\xi(H(s)))$,
\item $(\forall\beta\in w^p)(\alpha(H(\beta,0))\in a^q\ \Leftrightarrow\ \beta
\in a^p)$, 
\item $(\forall s\in u^p)(f^p_s=f^q_{H(s)}\comp H)$,
\item $(\forall s_0,s_1\in u^p)(\alpha(s_0)\leq\alpha(s_1)\ \Rightarrow\
y_{s_0,s_1}=\{s\in u^p:H(s)\in y^q_{H(s_0),H(s_1)}\})$.
\end{enumerate}
 \end{definition}
 
\begin{proposition}
\label{4.2}
Assume $\kappa^{<\kappa}=\kappa$. Then $\bPk$ is a $\kappa$--complete
$\kappa^+$--cc forcing notion of size $\kappa^+$.
\end{proposition}

\Proof It should be clear that $\bPk$ is $\kappa$--complete and $|\bPk|=
\kappa^+$. Moreover, there is $\kappa$ many isomorphism types of conditions in
$\bPk$ (and a condition in $\bPk$ is determined by its isomorphism type and
the set $w^p$). Now, to show the chain condition assume that $\cA\subseteq\bPk$
is of size $\kappa^+$. Applying $\Delta$--lemma choose pairwise isomorphic
conditions $p_0,p_1,p_2\in\cA$ such that $\{w^{p_0},w^{p_1},w^{p_2}\}$ forms a
$\Delta$--system with heart $w^*$ and such that for $i<j<3$
\[\sup(w^*)<\min(w^{p_i}\setminus w^*)\leq\sup(w^{p_i})<\min(w^{p_j}\setminus
w^*)\] 
(remember $\kappa^{<\kappa}=\kappa$). For $i,j<3$ let $H_{i,j}:u^{p_i}
\longrightarrow u^{p_j}$ be the isomorphism from $p_i$ to $p_j$. We are going
to define a condition $q\in\bPk$ which will be an upper bound to $p_1,p_2$
(note: not $p_0$!). To this end we first let
\[w^q=w^{p_0}\cup w^{p_1}\cup w^{p_2},\quad u^q=u^{p_0}\cup u^{p_1}\cup
u^{p_2},\quad a^q=a^{p_1}\cup a^{p_2}.\] 
To define functions $f^q_s$ we use the approach which can be described as ``put
zero whenever possible''. Thus we let: 
\begin{itemize}
\item if $s\in u^{p_1}\setminus u^{p_0}$ then $f^q_s={\bf 0}_{u^{p_0}}\cup
f^{p_1}_s\cup {\bf 0}_{u^{p_2}}$, 
\item if $s\in u^{p_2}\setminus u^{p_0}$ then $f^q_s={\bf 0}_{u^p_0}\cup {\bf
0}_{u^p_1}\cup f^{p_2}_s$,
\item if $s\in u^{p_0}$ then $f^q_s=f^{p_0}_s\cup f^{p_1}_{H_{0,1}(s)}\cup
f^{p_2}_{H_{0,2}(s)}$. 
\end{itemize}
It should be clear that the functions $f^q_s$ are well defined. Now we are
going to define the sets $y^q_{s_0,s_1}$ for distinct $s_0,s_1\in u^q$ such
that $\alpha(s_0)\leq\alpha(s_1)$. It is done by cases considering all
possible configurations. Thus we put:   
\begin{itemize}
\item if $s_0,s_1\in u^{p_i}$, $i<3$ then $y^q_{s_0,s_1}=y^{p_i}_{s_0,s_1}$,
\item if $s_0\in u^{p_1}\setminus u^{p_0}$, $s_1\in u^{p_2}\setminus u^{p_0}$
then $y^q_{s_0,s_1}=\{H_{2,0}(s_1)\}$,
\item if $s_0\in u^{p_0}$, $s_1\in u^{p_i}$, $i\in\{1,2\}$ then
\[y^q_{s_0,s_1}=\left\{\begin{array}{lll}
\emptyset & \mbox{if} & H_{i,0}(s_1)=s_0,\\
\{H_{i,0}(s_1)\} & \mbox{if} & \alpha(H_{i,0}(s_1))<\alpha(s_0),\\
y^{p_0}_{s_0,H_{i,0}(s_1)} & \mbox{if} & \alpha(s_0)\leq\alpha(H_{i,0}(s_1)),\
s_0\neq H_{i_,0}(s_1).\\
		       \end{array}\right.\]
\end{itemize}
We claim that 
\[q=\left\langle w^q,u^q,a^q,\langle f^q_s:s\in u^q\rangle,\langle y^q_{s_0,
s_1}:s_0,s_1\in u^q,\ s_0\neq s_1,\ \alpha(s_0)\leq\alpha(s_1)\rangle\right
\rangle\] 
is a condition in $\bPk$ and for this we have to check the demands of
\ref{4.1}. Clauses (a) and (b) should be obvious. To check \ref{4.1}(c) note
that $a^q\cap w^{p_0}=a^q\cap w^*$ and therefore there are no problems when
$\alpha\in a^q\cap w^{p_0}$. If $\alpha\in a^q\cap (w^{p_1}\setminus w^{p_0})$
and $\alpha<\beta\in a^q\cap (w^{p_2}\setminus w^{p_0})$ then $f^q_{\alpha,0}
(\beta,0)=f^q_{\alpha,1}(\beta,0)=0$. In all other instances we use the clause
(c) of \ref{4.1} for $p_1,p_2$.  

Now we have to verify the demand \ref{4.1}(d). Suppose that $s_0,s_1$ are
distinct members of $u^q$ and $\alpha(s_0)\leq\alpha(s_1)$. If $s_0,s_1\in
u^{p_i}$ for some $i<3$ then easily the set $y^q_{s_0,s_1}$ has the required
property. So suppose now that $s_0\in u^{p_1}\setminus u^{p_0}$, $s_1\in
u^{p_2}\setminus u^{p_0}$ (so then $f_{s_0}^q(s_1)=0$) and let $t\in u^q$ be
such that $f^q_t(s_1)=1=f^q_t(s_0)$. Then necessarily $t\in u^{p_0}$ and
$f^q_t(H_{2,0}(s_1))=f^q_t(s_1)=1$, so we are done in this case. Finally, let
us assume that $s_0\in u^{p_0}$ and $s_1\in u^{p_i}$, $0<i<3$. Note that if
$f^q_t(s_0)=1$ then $t\in u^{p_0}$. Now, if $H_{i,0}(s_1)=s_0$ then
$f^q_t(s_0)=f^q_t(s_1)$ for every $t\in u^{p_0}$ and there are no problems
(i.e.~no $f^q_t$ has to be taken care of). If $\alpha(H_{i,0}(s_1))<\alpha(
s_0)$ then the set $y^q_{s_0,s_1}=\{H_{i,0}(s_1)\}$ will work as for every
$t\in u^{p_0}$ we have $f^q_t(H_{i,0}(s_1))=f^q_t(s_1)$ (and $f^q_{s_0}(s_1)=
0$). For the same reason the set $y^q_{s_0,s_1}$ has the required property in
the remaining case too. 

Checking that the condition $q$ is stronger than both $p_1$ and $p_2$ is
straightforward (note: we do not claim that $q$ is stronger than $p_0$). \QED 

\begin{lemma}
\label{4.4}
If $p\in\bPk$, $t\in\kappa^+\times\kappa$ then there is $q\in\bPk$ such that
$p\leq q$ and $t\in u^q$.
\end{lemma}

\Proof Suppose $t=(\alpha,\xi)\in(\kappa^+\times\kappa)\setminus u^p$. Put
$w^q=w^p\cup\{\alpha\}$, $a^q=a^p$ and $u^q=u^p\cup\{(\alpha,0),(\alpha,1),
(\alpha,\xi)\}$. For $s\in u^p$ let $f^q_s=f^p_s\cup {\bf 0}_{u^q\setminus
u^p}$ and for $s\in u^q\setminus u^p$ let $f^q_s$ be such that $f^q_s(s)=1$
and $f^q_s\rest u^q\setminus\{s\}\equiv 0$. Finally, for distinct $s_0,s_1\in
u^q$ such that $\alpha(s_0)\leq\alpha(s_1)$ let 
\[y^q_{s_0,s_1}=\left\{\begin{array}{ll}
y^p_{s_0,s_1} & \mbox{if } s_0,s_1\in u^p,\\
\emptyset & \mbox{otherwise.}
		       \end{array}\right. \]
Check that $q=\langle w^q,u^q,a^q,\langle f^q_s:s\in u^q\rangle, \langle
y^q_{s_0,s_1}:s_0,s_1\in u^q\rangle\rangle\in \bPk$ is as required. \QED 
\medskip

For $p\in\bPk$ let $\bB_p$ be the algebra $\bB_{(u^p,F^p)}$ (see
\ref{1.1}(1)), where $F^p=\{f^p_s:s\in u^p\}\cup\{{\bf 0}_{u^p}\}$, and let
$\dot{\bB}_*$ be a $\bPk$--name such that  
\[\forces_{\bPk}\mbox{`` }\dot{\bB}_*=\bigcup\{\bB_p:p\in\Gamma_{\bPk}\}
\mbox{ ''.}\]
Furthermore, for $s\in\kappa^+\times\kappa$ let $\dot{f}_s$ be a $\bPk$--name
such that 
\[\forces_{\bPk}\mbox{`` }\dot{f}_s=\bigcup\{f^p_s:p\in\Gamma_{\bPk}\ \&\
s\in u^p\}\mbox{ ''.}\] 

\begin{proposition}
\label{4.5}
Assume $\kappa^{<\kappa}=\kappa$. Then in $\V^{\bPk}$:
\begin{enumerate}
\item $\dot{\bB}_*$ is the algebra $\bB_{(W,\dot{F})}$, where $W=\kappa^+
\times\kappa$ and $\dot{F}=\{\dot{f}_s: s\in\kappa^+\times\kappa\}\cup\{{\bf
0}_{\kappa^+\times\kappa}\}$, 
\item the algebra $\dot{\bB}_*$ is superatomic,
\item if $s\in\kappa^+\times\kappa$ and $b\in\dot{\bB}_*$ then there are
finite $v_0\subseteq v_1\subseteq \alpha(s)\times\kappa$ such that either
$x_s\wedge b$ or $x_s\wedge (-b)$ equals to 
\[\bigvee\{x_t\wedge\bigwedge_{\scriptstyle t'\in v_1\atop \scriptstyle
\alpha(t')<\alpha(t)} (-x_{t'}):\ \ t\in v_0\},\]
\item the height of $\dot{\bB}_*$ is $\kappa^+$ and $\{x_{\alpha,\xi}:\xi\in
\kappa\}$ are representatives of atoms of rank $\alpha+1$, 
\item $\irr(\dot{\bB}_*)=\kappa^+$.
\end{enumerate}
\end{proposition}

\Proof 1)\ \ \ First note that if $p\leq q$ then $\bB_p$ is a subalgebra of
$\bB_q$. Next, remembering \ref{4.4}, conclude that
\[\forces_{\bP^*}\mbox{`` $\dot{\bB}_*$ is a Boolean algebra generated by
}\langle x_s:s\in\kappa^+\times\kappa\rangle\mbox{ ''.}\]
Clearly, by \ref{4.4}, $\forces$`` $\dot{f}_s:\kappa^+\times\kappa
\longrightarrow\{0,1\}$ '' and $p\forces$`` $\dot{f}_s\rest u^p=f^p_s$ '' (for
$s\in u^p$, $p\in\bPk$). So it should be clear that $\forces_{\bPk}\dot{\bB}_*
=\bB_{(W,\dot{F})}$, where $W=\kappa^+\times\kappa$ and $\dot{F}=\{\dot{f}_s:
s\in \kappa^+\times\kappa\}\cup\{{\bf 0}_{\kappa^+\times\kappa}\}$.
\medskip

\noindent 2)\ \ \ It follows from \ref{4.1}(b) that for each $s\in\kappa^+
\times\kappa$ 
\[\forces_{\bPk}\mbox{`` }\dot{f}_s(s)=1\quad\mbox{and}\quad(\forall t\in
\kappa^+\times\kappa)(\alpha(t)\leq\alpha(s)\ \&\ t\neq s\ \ \Rightarrow\ \
\dot{f}_s(t)=0)\mbox{ ''.}\]
Now work in $\V^{\bPk}$. Let $\dot{J}_\alpha$ be the ideal in $\dot{\bB}_*$
generated by $\{x_{\beta,\xi}:\beta<\alpha, \xi\in\kappa\}$ (for $\alpha\leq
\kappa^+$; if $\alpha=0$ then $\dot{J}_\alpha=\{0\}$). It follows from the
previous remark that $x_{\alpha,\xi}\notin\dot{J}_\alpha$ (for all $\xi\in
\kappa$; remember \ref{1.2}). 

Suppose that $s_0,s_1$ are distinct, $\alpha(s_0)=\alpha(s_1)=\alpha<\kappa^+$
and suppose that $t\in\kappa^+\times\kappa$ is such that $\dot{f}_t(s_0)=
\dot{f}_t(s_1)=1$. Let $p\in\Gamma_{\bP^*}$ be such that $t,s_0,s_1\in u^p$. 
It follows from \ref{4.1}(d) that there is $s\in y^p_{s_0,s_1}$ such that
$f^p_t(s)=1$. Hence (applying \ref{1.2}) we may conclude that  
\[\dot{\bB}_*\models x_{s_0}\wedge x_{s_1}\leq\bigvee\{x_s: s\in y^p_{s_0,s_1}
\},\]
and therefore $x_{s_0}\wedge x_{s_1}\in\dot{J}_\alpha$. 

Now suppose that $s_0,s_1\in\kappa^+\times\kappa$ are such that $\alpha(s_0)<
\alpha(s_1)$ and let $p\in\Gamma_{\bPk}$ be such that $s_0,s_1\in u^p$. If
$f^p_{s_0}(s_1)=0$ then, by similar considerations as above, we have $x_{s_0}
\wedge x_{s_1}\in\dot{J}_\alpha$. Similarly, if $f^p_{s_0}(s_1)=1$ then $x_{
s_0}\wedge (-x_{s_1})\in\dot{J}_\alpha$. Hence we conclude that $x_{s_0}/
\dot{J}_\alpha$ is an atom in $\dot{\bB}_*/\dot{J}_\alpha$. 

Finally, note that the ideal $\dot{J}_{\kappa^+}$ is maximal (as $\{x_s:s\in
\kappa^+\times\kappa\}$ are generators of the algebra $\dot{\bB}_*$) and hence
the algebra $\dot{\bB}$ is superatomic.   
\medskip

\noindent 3)\ \ \ For $\alpha\leq\kappa^+$ let $\dot{J}_\alpha$ be the ideal
of $\dot{\bB}_*$ defined as above. Note that if $a\in\dot{J}_\alpha\setminus
\{0\}$ then there is a finite set $v\subseteq\alpha\times\kappa$ such that 
\[a\leq\bigvee_{t\in v}x_t\quad\mbox{ and }\quad (\forall t\in v)(x_t\wedge
a\notin\dot{J}_{\alpha(t)}).\]
A set $v$ with these properties will be called {\em a good $\alpha$--cover for
$a$}.

We know already that $x_s/\dot{J}_{\alpha(s)}$ is an atom in $\dot{\bB}_*/
\dot{J}_{\alpha(s)}$ and therefore either $x_s\wedge b\in\dot{J}_{\alpha(s)}$
or $x_s\wedge (-b)\in\dot{J}_{\alpha(s)}$. We may assume that the first takes
place. Applying repeatedly the previous remark find a finite set $v_1\subseteq
\alpha(s)\times\kappa$ such that for every $t\in v_1\cup\{s\}$:
\begin{enumerate}
\item if $x_t\wedge (x_s\wedge b)\in \dot{J}_{\alpha(t)}\setminus\{0\}$ then
there is a good $\alpha(t)$--cover $v\subseteq v_1$ for $x_t\wedge (x_s\wedge
b)$,
\item if $x_t\wedge (-x_s\vee -b)\in \dot{J}_{\alpha(t)}\setminus\{0\}$ then
there is a good $\alpha(t)$--cover $v\subseteq v_1$ for $x_t\wedge (-x_s\vee
-b)$.
\end{enumerate}
Now let $v_0=\{t\in v_1: x_t\wedge (x_s\wedge b)\notin \dot{J}_{\alpha(t)}\}$
and check that 
\[x_s\wedge b=\bigvee\{x_t\wedge\bigwedge_{\scriptstyle t'\in v_1\atop
\scriptstyle \alpha(t')<\alpha(t)} (-x_{t'}):\ \ t\in v_0\},\]
as required.
\medskip

\noindent 4)\ \ \ Almost everything what we need for this conclusion was done
in clause 2) above except that we have to check that, for each $\alpha<
\kappa^+$, $\{x_{\alpha,\xi}/\dot{J}_\alpha: \xi<\kappa\}$ lists {\em all}
atoms of the algebra $\dot{\bB}_*/\dot{J}_\alpha$. So suppose that $b/
\dot{J}_{\alpha}$ is an atom in $\dot{\bB}_*/\dot{J}_\alpha$. We may assume
that $b=\bigwedge\limits_{t\in w} x_t\wedge\bigwedge\limits_{t\in u} (-x_t)$
and that $\alpha(t)>\alpha$ for $t\in w$ (otherwise either $b\in
\dot{J}_\alpha$ or $b/\dot{J}_\alpha=x_s/\dot{J}_\alpha$ for some $s$ with
$\alpha(s)=\alpha$).  
   
Suppose that $w=\emptyset$. Let $p\in\bPk$. We may find a condition $q\geq p$
such that $u\subseteq u^q$ and then take $t\in (\{\alpha\}\times\kappa)
\setminus u^q$. Exactly as in the proof of \ref{4.4} we define a condition
$r\in \bPk$ stronger than $q$ and such that $t\in u^r$. Note that for this
condition we have $r\forces x_t\leq b$ and we easily finish. 

Let $s\in w$ (so $\alpha(s)>\alpha$) and $b^*=\bigwedge\limits_{t\in w
\setminus \{s\}}x_t\wedge\bigwedge\limits_{t\in u} (-x_t)$ (so $b=b^*\wedge
x_s$). It follows from the third clause that we find finite sets $v_0\subseteq
v_1\subseteq\alpha(s)\times\kappa$ such that 
\[c\stackrel{\rm def}{=}\bigvee\{x_t\wedge\bigwedge_{\scriptstyle t'\in v_1
\atop \scriptstyle\alpha(t')<\alpha(t)} (-x_{t'}):\ \ t\in v_0\}\in
\{x_s\wedge b^*,x_s\wedge (-b^*)\}.\]
If $c=x_s\wedge (-b^*)$ then we repeat arguments similar to those from the
previous paragraph but with a modified version of \ref{4.4}: defining the
condition $r$ with the property that $t\in u^r$, we use the function $f^q_s
\cup\{(t,1)\}$ as $f^r_t$ (check that no changes are needed in the definition
of $y^r_{s_0,s_1}$). Then easily $r\forces x_t\leq x_s\wedge (-c)$.  Finally,
if $c=x_s\wedge b^*$ then we take $s'\in v_0$ such that $\alpha(s')$ is
maximal possible. If $\alpha(s')>\alpha$ then similarly as in the previous
case we find a condition $r$ which forces that $x_t\leq x_s\wedge b^*=b$, if
$\alpha(s')\leq\alpha$ it is even easier. In all cases we easily finish
finding an element $x_{\alpha,\zeta}$ which is $\dot{J}_\alpha$--smaller than
$b$. 
\medskip

\noindent 5)\ \ \ Look at the demand \ref{4.1}(c): it means that if $\alpha,
\beta\in\dot{a}\stackrel{\rm def}{=}\bigcup\{a^p:p\in\Gamma_{\bPk}\}$ are
distinct then $\dot{f}_{\alpha,0}(\beta,0)=\dot{f}_{\alpha,1}(\beta,0)$. As
$\dot{f}_{\alpha,0}(\alpha,0)=1$, $\dot{f}_{\alpha,1}(\alpha,0)=0$ we conclude
that $\dot{f}_{\alpha,0},\dot{f}_{\alpha,1}$ determine homomorphisms from
$\dot{\bB}_*$ to $\{0,1\}$ witnessing $x_{\alpha,0}\notin\langle x_{\beta,0}:
\beta\in\dot{a}\setminus\{\alpha\}\rangle_{\dot{\bB}_*}$. Since clearly
$\forces |\dot{a}|=\kappa^+$ the proof is finished. \QED 

\begin{proposition}
\label{4.6}
Assume $\kappa^{<\kappa}=\kappa$. Then
\[\forces_{\bPk} \inc(\dot{\bB}_*)=s(\dot{\bB}_*)=\kappa.\]
\end{proposition}

\Proof Suppose that $\langle\dot{b}_\alpha:\alpha<\kappa^+\rangle$ is a
$\bPk$--name for a $\kappa^+$--sequence of elements of $\dot{\bB}_*$ and 
\[p\forces_{\bPk}\mbox{`` }\langle\dot{b}_\alpha:\alpha<\kappa^+\rangle\mbox{
are pairwise incomparable ''.}\]
Applying $\Delta$--lemma and ``standard cleaning'' choose pairwise isomorphic
conditions $p_0,p_1,p_2$ stronger than $p$, sets $v_1,v_2$, a Boolean
term $\tau$ and $\alpha_1<\alpha_2<\kappa^+$ such that  
\begin{itemize}
\item $\{w^{p_0},w^{p_1},w^{p_2}\}$ forms a $\Delta$--system with heart $w^*$,
\item $\sup(w^*)<\min(w^{p_i}\setminus w^*)\leq\sup(w^{p_i})<\min(w^{p_j}
\setminus w^*)$ for $i<j<3$,
\item $v_i\in [u^{p_i}]^{\textstyle<\!\omega}$ for $i=1,2$,
\item if $H_{i,j}$ is the isomorphism from $p_i$ to $p_j$ then $v_2=H_{2,1}[ 
v_1]$,
\item $p_i\forces$`` $\dot{b}_{\alpha_i}=\tau(x_s:s\in v_i)$ '' for $i=1,2$.
\end{itemize}
Considering two cases, we are going to define a condition $r$ stronger than
$p_1,p_2$. The condition $r$ will be defined in a similar manner as the
condition $q$ in the proof of \ref{4.2}.  
\medskip

\noindent {\sc Case A}:\qquad $\{0,1\}\models\tau(0:t\in v_1)=0$.\\
First choose $s^*\in u^{p_2}\setminus u^{p_0}$ such that
\begin{quotation}
\noindent if there is $s\in u^{p_2}\setminus u^{p_0}$ with the property that
\[\{0,1\}\models\tau(f^{p_2}_s(t):t\in v_2)=1\]
then $s^*$ is like that.
\end{quotation}
Now we proceed as in \ref{4.2} using $f^{p_2}_{s^*}$ instead of ${\bf
0}_{u^p_2}$. So we let    
\[w^r=w^{p_0}\cup w^{p_1}\cup w^{p_2},\quad u^r=u^{p_0}\cup u^{p_1}\cup
u^{p_2},\quad a^r=a^{p_1}\cup a^{p_2},\]
and we define functions $f^r_s$ as follows:
\begin{itemize}
\item if $s\in u^{p_0}$ then $f^r_s=f^{p_0}_s\cup f^{p_1}_{H_{0,1}(s)}\cup
f^{p_2}_{H_{0,2}(s)}$,
\item if $s\in u^{p_1}\setminus u^{p_0}$ then $f^r_s={\bf 0}_{u^{p_0}}\cup
f^{p_1}_s\cup f^{p_2}_{s^*}$, 
\item if $s\in u^{p_2}\setminus u^{p_0}$ then $f^r_s={\bf 0}_{u^{p_0}}\cup
{\bf 0}_{u^{p_1}}\cup f^{p_2}_s$
\end{itemize}
(check that the functions $f^r_s$ are well defined). Next, for distinct $s_0,
s_1\in u^r$ such that $\alpha(s_0)\leq\alpha(s_1)$, we define the sets $y^r_{
s_0,s_1}$:
\begin{itemize}
\item if $s_0,s_1\in u^{p_i}$, $i<3$ then $y^r_{s_0,s_1}=y^{p_i}_{s_0,s_1}$,
\item if $s_0\in u^{p_1}\setminus u^{p_0}$, $s_1\in u^{p_2}\setminus u^{p_0}$
then $y^r_{s_0,s_1}=\{H_{1,0}(s_0)\}$,
\item if $s_0\in u^{p_0}$, $s_1\in u^{p_i}$, $i\in\{1,2\}$ then
\[y^r_{s_0,s_1}=\left\{\begin{array}{lll}
\emptyset & \mbox{if} & H_{i,0}(s_1)=s_0,\\
\{H_{i,0}(s_1)\} & \mbox{if} & \alpha(H_{i,0}(s_1))<\alpha(s_0),\\
y^{p_0}_{s_0,H_{i,0}(s_1)} & \mbox{if} & \alpha(s_0)\leq\alpha(H_{i,0}(s_1)),\
s_0\neq H_{i_,0}(s_1).\\
		       \end{array}\right.\]
\end{itemize}
Exactly as in \ref{4.2} one checks that  
\[r=\left\langle w^r,u^r,a^r,\langle f^r_s:s\in u^r\rangle,\langle y^r_{s_0,
s_1}:s_0,s_1\in u^r,\ s_0\neq s_1,\ \alpha(s_0)\leq\alpha(s_1)\rangle\right
\rangle\] 
is a condition in $\bPk$ stronger than both $p_1$ and $p_2$. Moreover, it
follows from the definition of $f^r_s$'s that 
\[\bB_r\models\tau(x_t:t\in v_1)\leq\tau(x_t:t\in v_2)\]
(see \ref{1.2}). Consequently $r\forces\dot{b}_{\alpha_1}\leq\dot{b}_{
\alpha_2}$, a contradiction. 
\medskip

\noindent {\sc Case B}:\qquad $\{0,1\}\models\tau(0:t\in v_1)=1$.\\
Define $r$ almost exactly like in Case A, except that when choosing
$s^*\in u^{p_2}\setminus u^{p_0}$ ask if there is $s\in u^{p_2}\setminus
u^{p_0}$ such that  
\[\{0,1\}\models\tau(f^{p_2}_s(t):t\in v_2)=0\]
(and if so then $s^*$ has this property). Continue like before getting a
condition $r$ stronger than $p_1,p_2$ and such that
\[\bB_r\models\tau(x_t:t\in v_1)\geq\tau(x_t:t\in v_2)\]
and therefore $r\forces\dot{b}_{\alpha_1}\geq\dot{b}_{\alpha_2}$, a
contradiction finishing the proof. \QED  

\begin{theorem}
\label{4.7}
Assume $\kappa^{<\kappa}=\kappa$. Then
\[\forces_{\bPk} \Id(\dot{\bB})=2^{\kappa}=(2^{\kappa})^{\V}.\]
\end{theorem}

\Proof Let $\cK$ be the a family of all pairs $(p,\tau)$ such that $p\in\bPk$
and $\tau=\tau(x_s:s\in v)$ is a Boolean term, $v\subseteq u^p$. For each
ordinal $\alpha<\kappa^+$ we define a relation $E^-_\alpha$ on $\cK$ as
follows: 
\medskip

$(p_0,\tau_0)\; E^-_\alpha\;(p_1,\tau_1)$\quad if and only if
\begin{enumerate}
\item[(i)]\ \ \ the conditions $p_0,p_1$ are isomorphic,
\item[(ii)]\ \  $w^{p_0}\cap\alpha=w^{p_1}\cap\alpha$,
\item[(iii)]\   if $H:u^{p_0}\longrightarrow u^{p_1}$ is the isomorphism from
$p_0$ to $p_1$ then $\tau_1=H(\tau_0)$ (i.e.~$\tau_0=\tau(x_s:s\in v)$,
$\tau_1=\tau(x_{H(s)}:s\in v)$).
\end{enumerate}
A relation $E_\alpha$ on $\cK$ is defined by
\medskip

$(p_0,\tau_0)\; E_\alpha\;(p_1,\tau_1)$\quad if and only if\quad
$(p_0,\tau_0)\; E^-_\alpha\;(p_1,\tau_1)$ and
\begin{enumerate}
\item[(iv)]\ \ \ if $\beta\in w^{p_0}$ then

$\beta-\sup(w^{p_0}\cap\beta)=H(\beta)-\sup(w^{p_1}\cap H(\beta))\ \mod\
\kappa$\qquad and

$\beta\geq\sup(w^{p_0}\cap\beta)+\kappa$ if and only if $H(\beta)\geq\sup(
w^{p_1}\cap H(\beta))+\kappa$.
\end{enumerate}

\begin{claim}
\label{4.7.1}
For each $\alpha<\kappa^+$, $E_\alpha$, $E^-_\alpha$ are equivalence relations
on $\cK$ with $\kappa$ many equivalence classes.
\end{claim}

\begin{claim}
\label{4.7.2}
Suppose that $\alpha<\kappa^+$, $(p_0,\tau_0)\; E_\alpha\; (p_1,\tau_1)$ and
$p_0\leq q_0$. Then there is $q_1\in\bPk$ such that $p_1\leq q_1$ and
$(q_0,\tau_0)\; E^-_\alpha\; (q_1,\tau_1)$.
\end{claim}

\begin{claim}
\label{4.7.3}
Suppose that $\dot{I}$ is a $\bPk$--name for an ideal in the algebra
$\dot{\bB}_*$ and let $\cK(\dot{I})=\{(p,\tau)\in\cK:p\forces\tau\in\dot{I}
\}$. Then there is $\alpha=\alpha(\dot{I})<\kappa^+$ such that
\[\cK(\dot{I})=\bigcup\{(p,\tau)/E_\alpha: (p,\tau)\in\cK(\dot{I})\}.\]
\end{claim}

\noindent{\em Proof of the claim:}\qquad Assume not. Then for each
$\alpha<\kappa^+$ we find $(p^\alpha_0,\tau^\alpha_0)\in\cK(\dot{I})$ and
$(p^\alpha_1,\tau^\alpha_1)\notin\cK(\dot{I})$ such that $(p^\alpha_0,
\tau^\alpha_0)\; E_\alpha\;(p^\alpha_1,\tau^\alpha_1)$. Take $q^\alpha_1\geq
p^\alpha_1$ such that $q^\alpha_1\forces\tau^\alpha_1\notin\dot{I}$ and use
\ref{4.7.2} to find $q^\alpha_0\geq p^\alpha_0$ such that $(q^\alpha_0,
\tau^\alpha_0)\; E^-_\alpha\;(q^\alpha_1,\tau^\alpha_1)$. Now use
$\Delta$--lemma and clause (i) of the definition of $E^-_\alpha$ to find
$\alpha_0<\alpha_1<\alpha_2<\alpha_3<\kappa^+$ such that letting $q_2=q^{
\alpha_2}_1$, $\tau_2=\tau^{\alpha_2}_1$ and $q_i=q^{\alpha_i}_0$, $\tau_i=
\tau^{\alpha_i}_0$ for $i\neq 2$ we have 
\begin{itemize}
\item the conditions $q_0,\ldots,q_3$ are pairwise isomorphic (and for $i,j<4$
let $H_{i,j}:u^{q_i}\longrightarrow u^{q_j}$ be the isomorphism from $q_i$ to
$q_j$), 
\item $\{w^{q_0}, w^{q_1}, w^{q_2},w^{q_3}\}$ forms a $\Delta$--system with
heart $w^*$, 
\item $\sup(w^*)<\min(w^{q_i}\setminus w^*)\leq\sup(w^{q_i}\setminus w^*)
<\min(w^{q_j}\setminus w^*)$ when $i<j<4$,
\item $\tau_i=H_{i,j}(\tau_j)$ (i.e.~we have the same term). 
\end{itemize}
Now we define a condition $q\in\bPk$ in a similar manner as in \ref{4.2},
\ref{4.6}. First we fix $s^*\in u^{q_3}\setminus u^{q_0}$ such that
\begin{quotation}
\noindent if there is $s\in u^{q_3}\setminus u^{q_0}$ with the property that
$f^{q_3}_s(\tau_3)=1$\\
then $s^*$ is like that.
\end{quotation}
We put
\[w^q=w^{q_0}\cup\ldots\cup w^{q_3},\quad u^q=u^{q_0}\cup\ldots\cup u^{q_3},
\quad a^q=a^{q_1}\cup a^{q_2}\cup a^{q_3},\]
and we define $f^q_s$ as follows:
\[f^q_s=\left\{\begin{array}{lll}
\bigcup\limits_{i<4}f^{q_i}_{H_{0,i}(s)} &\mbox{ if}& s\in u^{q_0},\\
{\bf 0}_{u^{q_0}}\cup f^{q_1}_s\cup f^{q_2}_{H_{3,2}(s^*)}\cup f^{q_3}_{s^*}
& \mbox{ if}& s\in u^{q_1}\setminus u^{q_0},\\
{\bf 0}_{u^{q_0}}\cup {\bf 0}_{u^{q_1}}\cup f^{q_2}_s\cup f^{q_3}_{s^*}
& \mbox{ if}& s\in u^{q_2}\setminus u^{q_0},\\
{\bf 0}_{u^{q_0}}\cup {\bf 0}_{u^{q_1}}\cup {\bf 0}_{u^{q_2}}\cup f^{q_3}_s
& \mbox{ if}& s\in u^{q_3}\setminus u^{q_0}.\\
	       \end{array}\right.\]
Finally, for distinct $s_0,s_1\in u^q$ such that $\alpha(s_0)\leq\alpha(s_1)$,
we define 
\[y^q_{s_0,s_1}=\left\{\begin{array}{lll}
y^{q_i}_{s_0,s_1} &\mbox{ if}& s_0,s_1\in u^{q_i},\ i<4,\\
\{H_{i,0}(s_0)\}  &\mbox{ if}& s_0\in u^{q_i}\setminus u^{q_0},\ s_1\in u^{
q_j}\setminus u^{q_0},\ 0<i<j<4,\\
\emptyset         &\mbox{ if}& s_0\in u^{q_0},\ s_1\in u^{q_i},\ 0<i<4,\
H_{i,0}(s_1)=s_0,\\ 
\{H_{i,0}(s_1)\}  &\mbox{ if}& s_0\in u^{q_0},\ s_1\in u^{q_i},\ 0<i<4,\
\alpha(H_{i,0}(s_1))<\alpha(s_0),\\
y^{q_0}_{s_0,H_{i,0}(s_1)} & &\mbox{otherwise}.
		       \end{array}\right.\]
It should be a routine to check that this defines a condition $q\in\bPk$
stronger than $q_1,q_2,q_3$ and that (by \ref{1.2}) $\bB_q\models\tau_2\leq
\tau_1\vee \tau_3$ (remember that the terms are isomorphic). But this means
that   
\[q\forces_{\bPk}\mbox{`` }\tau^{\alpha_2}_1\leq \tau^{\alpha_1}_0\vee
\tau^{\alpha_3}_0\quad\&\quad \tau^{\alpha_2}_1\notin \dot{I}\quad\&\quad
\tau^{\alpha_1}_0,\tau^{\alpha_3}_0\in\dot{I}\mbox{ '',}\]
a contradiction finishing the proof of the claim.
\medskip

Now, using \ref{4.7.3}, we may easily finish: if $\dot{I}_0$, $\dot{I}_1$ are
$\bPk$--names for ideals in $\dot{\bB}_*$ such that $\cK(\dot{I}_0)=\cK(
\dot{I}_1)$ then $\forces\dot{I}_0=\dot{I}_1$. But \ref{4.7.3} says that
$\cK(\dot{I})$ is determined by $\alpha(\dot{I})$ and a family of equivalence
classes of $E_{\alpha(\dot{I})}$. So we have at most $\kappa^+\cdot 2^\kappa=
2^\kappa$ possibilities for $\cK(\dot{I})$. Finally note that $|\bPk|=
\kappa^+$ and $\bPk$ satisfies the $\kappa^+$--cc, so $\forces_{\bPk}2^\kappa
=(2^\kappa)^{\V}$. \QED

\begin{conclusion}
It is consistent that there is a superatomic Boolean algebra $\bB$ such that
\[s(\bB)=\inc(\bB)=\kappa,\quad\irr(\bB)=\Id(\bB)=\kappa^+\quad\mbox{and}\quad
\Sub(\bB)=2^{\kappa^+}.\qquad \QED\]
\end{conclusion}

\section{Modifications of $\bPk$}
In this section we modify the forcing notion $\bPk$ of \ref{4.1} and we get
two new models. The first model shows the consistency of ``there is a
superatomic Boolean algebra $\bB$ such that $\irr(\bB)<\inc(\bB)$'' answering
\cite[Problem 79]{M2}. Next we solve \cite[Problem 81]{M2} showing that
possibly there is a superatomic Boolean algebra $\bB$ with $\aut(\bB)<\tig(
\bB)$. 

\begin{definition}
\label{5.1}
Let $\kappa$ be a cardinal. A forcing notion $\bPz$ is defined like $\bPk$ of
\ref{4.1} but the demand \ref{4.1}(c) is replaced by: 
\begin{enumerate}
\item[(c$^0$)] \quad if $\alpha<\beta$, $\alpha,\beta\in a^p$ then\quad 
$(\exists s\in u^p)(f^p_s(\alpha,0)=1\ \&\ f^p_s(\beta,0)=0)$.
\end{enumerate}
\end{definition}

Naturally we have a variant of definition \ref{isomorph} of isomorphic
conditions for the forcing notion $\bPz$ (with no changes) and similarly as
for the case of $\bPk$ we define algebras $\bB_p$ (for $p\in\bPz$) and
$\bPz$--names $\dot{\bB}^0_*$, $\dot{f}^0_s$ (for $s\in\kappa^+\times\kappa$). 

\begin{proposition}
\label{5.2}
Assume $\kappa^{<\kappa}=\kappa$. Then $\bPk$ is a $\kappa$--complete
$\kappa^+$--cc forcing notion of size $\kappa^+$.
\end{proposition}

\Proof Repeat the proof of \ref{4.2} (with no changes). \QED

\begin{proposition}
\label{5.3}
Assume $\kappa^{<\kappa}=\kappa$. Then in $\V^{\bPk}$:
\begin{enumerate}
\item $\dot{\bB}^0_*$ is the algebra $\bB_{(W,\dot{F})}$, where $W=\kappa^+
\times\kappa$ and $\dot{F}=\{\dot{f}^0_s: s\in\kappa^+\times\kappa\}\cup\{{\bf
0}_{\kappa^+\times\kappa}\}$, 
\item the algebra $\dot{\bB}^0_*$ is superatomic (of height $\kappa^+$) and
$\{x_{\alpha,\xi}:\xi\in\kappa\}$ are representatives of atoms of rank
$\alpha+1$,  
\item $\inc(\dot{\bB}^0_*)=\kappa^+$.
\end{enumerate}
\end{proposition}

\Proof The proofs of the first two clauses are repetitions of that of
\ref{4.5}(1--4) (so we have the respective version of \ref{4.5}(3) too).

To show the third clause let $\dot{a}\stackrel{\rm def}{=}\bigcup\{a^p:p\in
\Gamma_{\bPz}\}$. It should be clear that $\forces |\dot{a}|=\kappa^+$. Note
that if $\alpha,\beta\in a^p$, $\alpha<\beta$ then, by \ref{5.1}(c$^0$),
$\bB_p\models x_{\alpha,0}\not\leq x_{\beta,0}$ and by the respective variant
of \ref{4.1}(b) we have $\bB_p\models x_{\beta,0}\not\leq x_{\alpha,0}$. 
Consequently the sequence $\langle x_{\alpha,0}:\alpha\in\dot{a}\rangle$
witnesses $\inc(\dot{\bB}^0_*)=\kappa^+$. \QED 

\begin{proposition}
\label{5.4}
Assume $\kappa^{<\kappa}=\kappa$. Then\quad $\forces_{\bPz}\irr^+_3(\dot{
\bB}^0_*)=\kappa^+$.
\end{proposition}

\Proof Let $\langle\dot{b}_\alpha:\alpha<\kappa^+\rangle$ be a $\bPz$--name
for a $\kappa^+$--sequence of elements of $\dot{\bB}^0_*$ and let $p\in\bPz$. 
Find pairwise isomorphic conditions $p_i$, sets $v_i$, ordinals $\alpha_i$
(for $i<7$) and a Boolean term $\tau$ such that  
\begin{itemize}
\item $p\leq p_0,\ldots,p_7$,\ \ $\alpha_0<\alpha_1<\ldots<\alpha_6<\kappa^+$,
\ \ $v_i\in [u^{p_i}]^{\textstyle <\!\omega}$ for $i<7$,
\item $\{w^{p_0},\ldots,w^{p_6}\}$ forms a $\Delta$--system with heart $w^*$,
\item $\sup(w^*)<\min(w^{p_i}\setminus w^*)\leq\sup(w^{p_i})<\min(w^{p_j}
\setminus w^*)$ for $i<j<7$,
\item if $H_{i,j}$ is the isomorphism from $p_i$ to $p_j$ then $v_j=H_{i,j}[ 
v_i]$ (for $i,j<7$),
\item $p_i\forces$`` $\dot{b}_{\alpha_i}=\tau(x_s:s\in v_i)$ '' for $i<7$.
\end{itemize}
Now we are going to define an upper bound $q$ to the conditions $p_3,\ldots,
p_6$. For this we let 
\[w^q=\bigcup_{i<7} w^{p_i},\quad u^q=\bigcup_{i<7} u^{p_i},\quad a^q=
\bigcup_{2<i<7} a^{p_i}\] 
and for $s\in w^q$ we define
\[f^q_s=\left\{\begin{array}{lll}
\bigcup\limits_{j<7}f^{p_j}_{H_{0,j}(s)} &\mbox{if}& s\in u^{p_0},\\
{\bf 0}_{u^{p_0}\cup u^{p_2}\cup u^{p_4}}\cup f^{p_1}_s\cup f^{p_3}_{H_{1,3}
(s)}\cup f^{p_5}_{H_{1,5}(s)}\cup f^{p_6}_{H_{1,6}(s)} &\mbox{if}& s\in
u^{p_1}\setminus u^{p_0},\\ 
{\bf 0}_{u^{p_0}\cup u^{p_1}\cup u^{p_5}}\cup f^{p_2}_s\cup f^{p_3}_{H_{2,3}
(s)}\cup f^{p_4}_{H_{2,4}(s)}\cup f^{p_6}_{H_{2,6}(s)} &\mbox{if}& s\in
u^{p_2}\setminus u^{p_0},\\ 
{\bf 0}_{u^{p_0}\cup u^{p_1}\cup u^{p_2}\cup u^{p_6}}\cup f^{p_3}_s\cup
f^{p_4}_{H_{3,4}(s)}\cup f^{p_5}_{H_{3,5}(s)} &\mbox{if}& s\in
u^{p_3}\setminus u^{p_0},\\ 
{\bf 0}_{u^q\setminus u^{p_i}}\cup f^{p_i}_s &\mbox{if}& s\in
u^{p_i}\setminus u^{p_0}, 3<i.\\ 
	       \end{array}\right.\]
Next, for distinct $s_0,s_1\in u^q$ such that $\alpha(s_0)\leq\alpha(s_1)$,
we define $y^q_{s_0,s_1}$ considering all possible configurations
separately. Thus we put:
\begin{itemize}
\item if $s_0,s_1\in u^{p_i}$, $i<7$ then $y^q_{s_0,s_1}=y^{p_i}_{s_0,s_1}$,
\item if $s_0\in u^{p_0}\setminus u^{p_1}$, $s_1\in u^{p_i}\setminus u^{p_0}$,
$0<i<7$ then 
\[y^q_{s_0,s_1}=\left\{\begin{array}{lll}
\emptyset &\mbox{if}& H_{i,0}(s_1)=s_0,\\
\{H_{i,0}(s_1)\} &\mbox{if}& \alpha(H_{i,0}(s_1))<\alpha(s_0),\\
y^{p_0}_{s_0,H_{i,0}(s_1)}& &\mbox{otherwise,}
		       \end{array}\right.\]
\item if $s_0\in u^{p_i}\setminus u^{p_0}$, $s_1\in u^{p_j}\setminus u^{p_0}$,
$0<i<j<7$ then 
\[y^q_{s_0,s_1}=\left\{\begin{array}{lll}
\{H_{i,k}(s_0):k<i\} &\mbox{if}& H_{j,i}(s_1)=s_0,\\
\{H_{i,k}(s_0):k<i\}\cup \{H_{j,i}(s_1)\} &\mbox{if}& \alpha(H_{j,i}(s_1))<
\alpha(s_0),\\
\{H_{i,k}(s_0):k<i\}\cup y^{p_i}_{s_0,H_{j,i}(s_1)}& &\mbox{otherwise.}
		       \end{array}\right.\]
\end{itemize}
It is not difficult to check that the above formulas define a condition
$q\in\bPz$ stronger than $p_3,p_4,p_5,p_6$ (just check all possible
cases). Moreover, applying \ref{1.2}, one sees that
\[\begin{array}{lll}
\bB_q\models &\tau(x_s:s\in v_3)=&(\tau(x_s:s\in v_4)\wedge\tau(x_s:s\in v_5))
\ \vee\\
& &(\tau(x_s:s\in v_4)\wedge\tau(x_s:s\in v_6))\ \vee\\
& &(\tau(x_s:s\in v_5)\wedge\tau(x_s:s\in v_6)).
  \end{array}\]
Hence
\[q\forces_{\bPz}\mbox{`` }\dot{b}_{\alpha_3}\in\langle\dot{b}_{\alpha_4},
\dot{b}_{\alpha_5}, \dot{b}_{\alpha_6}\rangle_{\dot{\bB}^0_*}\mbox{ ''},\]
finishing the proof. \QED

\begin{conclusion}
It is consistent that there is a superatomic Boolean algebra $\bB$ such that
$\inc(\bB)=\kappa^+$ and $\irr(\bB)=\kappa$. \QED
\end{conclusion}

For the next model we need a more serious modification of $\bPk$ involving a
change in the definition of the order. 

\begin{definition}
\label{5.5}
For an uncountable cardinal $\kappa$ we define a forcing notion $\bPj$ like
$\bPk$ of \ref{4.1} except that the clause \ref{4.1}(c) is replaced by
\begin{enumerate}
\item[(c$^1$)] \quad if $\alpha<\beta$, $\alpha,\beta\in a^p$ then\ \ 
$f^p_{\alpha,0}(\beta,0)=1$
\end{enumerate}
and we add the following requirement 
\begin{enumerate}
\item[(e)] \quad if $(1,\xi)\in u^p$ then the set $\{\zeta<\kappa: f^p_{0,
\zeta}(1,\xi)=1\}$ is infinite.
\end{enumerate}
Moreover, we change the definition of the order demanding additionally that,
if $p\leq q$,
\begin{enumerate}
\item[$(\alpha)$] if $(1,\xi)\in u^p$, $(0,\zeta)\in u^q\setminus u^p$ then
$f^q_{0,\zeta}(1,\xi)=0$, and 
\item[$(\beta)$]  if $(1,\xi)\in u^q\setminus u^p$ then the set $\{(0,\zeta)
\in u^p: f^p_{0,\zeta}(1,\xi)=1\}$ is finite.
\end{enumerate}
\end{definition}
Like before we have the respective variants of \ref{4.2}--\ref{4.5} for
$\bPj$ which we formulate below. The $\bPj$--names $\dot{\bB}^1_*$ and
$\dot{f}^1_s$ are defined like $\dot{\bB}_*$ and $\dot{f}_s$.

\begin{proposition}
\label{5.6}
Assume $\omega_0<\kappa=\kappa^{<\kappa}$. Then $\bPj$ is a $\kappa$--complete
$\kappa^+$--cc forcing notion.
\end{proposition}

\Proof Repeat the arguments of \ref{4.2} with the following small adjustments.
First note that we may assume $|w^*|>2$. Next, if $a^{p_2}\setminus w^*\neq
\emptyset$ then we let $\alpha=\min(a^{p_2}\setminus w^*)$ and defining
$f^q_s$ for $s\in u^{p_1}\setminus u^{p_0}$ we put $f^q_s= {\bf 0}_{u^{p_0}}
\cup f^{p_1}_s\cup f^{p_2}_{\alpha,0}$. (No other changes needed.) \QED

\begin{proposition}
\label{5.7}
Assume $\omega_0<\kappa=\kappa^{<\kappa}$. Then in $\V^{\bPj}$:
\begin{enumerate}
\item $\dot{\bB}^1_*$ is the algebra $\bB_{(W,\dot{F})}$, where $W=\kappa^+
\times\kappa$ and $\dot{F}=\{\dot{f}^1_s: s\in\kappa^+\times\kappa\}\cup\{{\bf
0}_{\kappa^+\times\kappa}\}$, 
\item the algebra $\dot{\bB}^1_*$ is superatomic,
\item if $s\in\kappa^+\times\kappa$ and $b\in\dot{\bB}^1_*$ then there are
finite $v_0\subseteq v_1\subseteq \alpha(s)\times\kappa$ such that either
$x_s\wedge b$ or $x_s\wedge (-b)$ equals to 
\[\bigvee\{x_t\wedge\bigwedge_{\scriptstyle t'\in v_1\atop \scriptstyle
\alpha(t')<\alpha(t)} (-x_{t'}):\ \ t\in v_0\},\]
\item the height of $\dot{\bB}^1_*$ is $\kappa^+$ and $\{x_{\alpha,\xi}:\xi\in
\kappa\}$ are representatives of atoms of rank $\alpha+1$, 
\item $\tig(\dot{\bB}^1_*)=\kappa^+$.
\end{enumerate}
\end{proposition}

\Proof (1)--(3) Repeat the arguments of \ref{4.5}(1--3) with no changes.
\medskip

\noindent (4) Like \ref{4.5}(4), but the cases $\alpha=0$ and $\alpha=1$ are
considered separately (for $\alpha>1$ no changes are required). 
\medskip

\noindent (5) Let $\dot{a}\stackrel{\rm def}{=}\bigcup\{a^p:p\in\Gamma_{\bPk}
\}$ and look at the sequence $\langle -x_{\alpha,0}:\alpha\in\dot{a}\rangle$. 
It easily follows from \ref{5.5}(c$^1$) that it is a free sequence (so it
witnesses $\tig(\dot{\bB}^1_*)=\kappa^+$). \QED

\begin{theorem}
\label{5.8}
Assume $\omega_0<\kappa=\kappa^{<\kappa}$. Then $\forces_{\bPj}$`` $\aut(
\dot{\bB}^1_*)=\kappa$ ''.
\end{theorem}

\Proof It follows from \ref{5.6} that, in $\V^{\bPj}$, $\kappa=\kappa^{<
\kappa}$. By \ref{5.7}(2,4) we have that each automorphism of $\dot{\bB}^1_*$
is determined by its values on atoms of $\dot{\bB}^1_*$ and $\{x_{0,\xi}:\xi<
\kappa\}$ is the list of the atoms of $\dot{\bB}^1_*$. Therefore it is enough
to show that in $\V^{\bPj}$:
\smallskip

if $\dot{h}:\dot{\bB}^1_*\longrightarrow\dot{\bB}^1_*$ is an automorphism
then $|\{\xi<\kappa: \dot{h}(x_{0,\xi})\neq x_{0,\xi}\}|<\kappa$.
\smallskip

\noindent So assume that $\dot{h}$ is a $\bPj$--name for an automorphism of
the algebra 
$\dot{\bB}^1_*$ and $p\in\bPj$ is such that $0,1\in w^p$. Now we consider three
cases. 
\medskip

\noindent {\sc Case A:}\quad for each $q\geq p$ there are $r\in\bPj$ and
distinct $\xi,\zeta<\kappa$ such that 
\[q\leq r,\quad (0,\xi),(0,\zeta)\in u^r\setminus u^q,\quad f^r_{0,\xi}\rest
u^q\equiv {\bf 0}\quad\mbox{and}\quad r\forces_{\bPj}\mbox{`` }\dot{h}(x_{0,
\xi})=x_{0,\zeta}\mbox{ ''}.\]
Construct inductively a sequence $\langle q_n,\xi_n,\zeta_n: n<\omega\rangle$
such that
\begin{itemize}
\item $q_n\in\bPj$,\quad $\xi_n,\zeta_n<\kappa$,\quad $p=q_0\leq q_1\leq q_2
\leq\ldots$, 
\item $(0,\xi_n),(0,\zeta_n)\in u^{q_{n+1}}\setminus u^{q_n}$\quad and\quad
$f^{q_{n+1}}_{0,\xi_n}\rest u^{q_n}\equiv {\bf 0}$,
\item $q_{n+1}\forces$`` $\dot{h}(x_{0,\xi_n})=x_{0,\zeta_n}$ ''.
\end{itemize}
Choose $\xi<\kappa$ such that $(1,\xi)\notin \bigcup\limits_{n<\omega}
u^{q_n}$. Now we are defining a condition $r\in\bPj$. First we put
\[w^r=\bigcup_{n<\omega} w^{q_n},\quad u^r=\{(1,\xi)\}\cup\bigcup_{n<\omega}
u^{q_n}\quad \mbox{and}\quad a^r=\bigcup_{n<\omega} a^{q_n}.\]
Next for $s\in u^q$ we put
\[f^r_s=\left\{\begin{array}{lll}
\{\langle (1,\xi),1\rangle\}\cup \bigcup\limits_{m>n} f^{q_m}_s&\mbox{if}&
s=(0,\xi_n),\ n\in\omega,\\
\{\langle (1,\xi),0\rangle\}\cup \bigcup\limits_{m>n} f^{q_m}_s&\mbox{if}&
s\in u^{q_n}\setminus\{(0,\xi_\ell):\ell\leq n\},\ n\in\omega,\\
{\bf 0}_{u^r\setminus\{s\}}\cup \{\langle s,1\rangle\}&\mbox{if}&s=(1,\xi).\\
	       \end{array}\right.\]
Furthermore, if $s_0,s_1\in u^r$ are distinct and such that $\alpha(s_0)\leq
\alpha(s_1)$ then we define $y^r_{s_0,s_1}$ as follows:
\medskip

\noindent{--} if $(1,\xi)\notin\{s_0,s_1\}$ then $y^r_{s_0,s_1}=y^{q_n}_{s_0,
s_1}$, where $n<\omega$ is such that $s_0,s_1\in u^{q_n}$, 

\noindent{--} if $s_0=(1,\xi)$, $s_1\in u^{q_n}$, $n<\omega$ then $y^r_{s_0,
s_1}=\{(0,\xi_m): m\leq n\}$, 

\noindent{--} if $s_1=(1,\xi)$, $s_0\in u^{q_n}$, $\alpha(s_0)=1$, $n<\omega$
then $y^r_{s_0,s_1}=\{(0,\xi_m): m\leq n\}$. 
\medskip

\noindent It is not difficult to check that the above formulas define a
condition $r\in\bPj$ stronger than all $q_n$ (verifying \ref{4.1}(d) remember
that $f^{q_{n+1}}_{0,\xi_n}\rest u^{q_n}\equiv {\bf 0}$). Note that $r\forces
(\forall n<\omega)(x_{0,\xi_n}\leq x_{1,\xi})$ and hence $r\forces (\forall
n<\omega)(x_{0,\zeta_n}\leq \dot{h}(x_{1,\xi}))$. Take a condition $r^*$
stronger than $r$ and such that for some $\zeta<\kappa$ we have $(1,\zeta)\in
u^{r^*}$ and $r^*\forces\dot{h}(x_{1,\xi})/\dot{J}_1=x_{1,\zeta}/\dot{J}_1$,
where $\dot{J}_1$ is the ideal of $\dot{\bB}^1_*$ generated by atoms (remember
\ref{5.7}(4)). Then for some $N$ we have $r^*\forces (\forall n\geq N)(x_{0,
\zeta_n}\leq x_{1,\zeta})$. Now look at the definition of the order in $\bPj$: 
by \ref{5.5}($\beta$) we have $(1,\zeta)\in u^r$. If $(1,\zeta)\in u^{q_n}$
for some $n<\omega$ then we get immediate contradiction with
\ref{5.5}($\alpha$), so the only possibility is that $\xi=\zeta$. But then
look at the definition of the functions $f^r_{0,\zeta_m}$ -- they all take
value 0 at $(1,\xi)$ so $r\forces x_{0,\zeta_n}\not\leq x_{1,\xi}$, a
contradiction.  Thus necessarily Case A does not hold.
\medskip

\noindent {\sc Case B:}\quad there are $p^*\geq p$ and $t\in u^{p^*}$
such that

for each $q\geq p^*$ there are $r\in\bPj$ and distinct $\xi,\zeta<\kappa$ with:

$q\leq r$, $(0,\xi),(0,\zeta)\in u^r\setminus u^q$, $r\forces_{\bPj}\mbox{``
}\dot{h}(x_{0,\xi})=x_{0,\zeta}\mbox{ ''}$, $f^r_{0,\xi}(t)=1$ and 

$(\forall s\in u^q)(\alpha(s)<\alpha(t)\ \Rightarrow\ f^r_{0,\xi}(s)=0)$.
\smallskip

\noindent First note that (by \ref{5.5}($\alpha$)) necessarily
$\alpha(t)>1$. Now apply the procedure of Case A with the following
modifications. Choosing $q_n,\xi_n,\zeta_n$ we demand that $q_0=p^*$,
$f^{q_{n+1}}_{0,\zeta_n}(t)=1$ and $(\forall s\in u^{q_n})(\alpha(s)<\alpha(t)
\ \Rightarrow\ f^{q_{n+1}}_{0,\xi_n}(t)=1)$. Next, defining the condition $r$
we declare that $f^r_{1,\xi}=\bigcup\limits_{n<\omega} f^{q_n}_t\cup\{\langle
(1,\xi),1\rangle\}$ and in the definition of $y^r_{s_0,s_1}$ we let 
\medskip

\noindent{--} if $s_0=(1,\xi)$ and either $s_1=t$ or $\alpha(s_1)<\alpha(t)$
then $y^r_{s_0,s_1}=\emptyset$,

\noindent{--} if $s_0=(1,\xi)$ and $\alpha(s_1)\geq\alpha(t)$, $s_1\neq t$
then $y^r_{s_0,s_1}=y^{q_n}_{t,s_1}$, where $n<\omega$ is such that $s_1\in
u^{q_n}$. 
\medskip

\noindent Continuing as in the Case A we get a contradiction.
\medskip

\noindent {\sc Case C:}\quad neither Case A nor Case B hold.
\smallskip

\noindent Let $q_0\geq p$ witness that Case A fails. So for each $r\geq q_0$
and distinct $\xi,\zeta<\kappa$ such that $(0,\xi),(0,\zeta)\in u^r\setminus
u^{q_0}$ 

if $r\forces\dot{h}(x_{0,\xi})=x_{0,\zeta}$ then $(\exists t\in u^{q_0})(
f^r_{0,\xi}(t)=1)$. 

\noindent Now, since Case B fails and $\bPj$ is $\kappa$--complete (and
$\kappa>\omega$) we may build a condition $q_1\geq q_0$ such that
\begin{quotation}
\noindent if $t\in u^{q_1}$, $r\geq q_1$, $(0,\xi),(0,\zeta)\in u^r\setminus
u^{q_1}$, $r\forces \dot{h}(x_{0,\xi})=x_{0,\zeta}$, $f^r_{0,\xi}(t)=1$ and
$(\forall s\in u^{q_1})(\alpha(s)<\alpha(t)\ \Rightarrow\ f^r_{0,\xi}(s)=0)$

\noindent then $\xi=\zeta$.
\end{quotation}
But then clearly 
\[q_1\forces_{\bPj}\mbox{`` }(\forall\xi<\kappa)(\dot{h}(x_{0,\xi})\neq x_{0,
\xi}\quad\Rightarrow\quad (0,\xi)\in u^{q_1})\mbox{ '',}\]
finishing the proof. \QED

\begin{conclusion}
It is consistent that there is a superatomic Boolean algebra $\bB$ such that
$\tig(\bB)=\kappa^+$ and $\aut(\bB)=\kappa$. \QED
\end{conclusion}

\section{When tightness is singular}
In this section we will show that, consistently, there is a Boolean algebra
with tightness $\lambda$ and such that there is an ultrafilter with this
tightness but there is no free sequence of length $\lambda$ and no homomorphic
image of the algebra has depth $\lambda$. This gives partial answers to
\cite[Problems 13, 41]{M2}. Next we show some bounds on possible consistency
results here showing that sometimes we may find quotients with depth equal to
the tightness of the original algebra.

Let us recall that a sequence $\langle b_\alpha: \alpha<\xi\rangle$ of
elements of a Boolean algebra $\bB$ is (algebraically) free if for each finite
sets $F,G\subseteq\xi$ such that $\max(F)<\min(G)$ we have 
\[\bB\models\bigwedge\limits_{\alpha\in F} b_\alpha\wedge\bigwedge\limits_{
\alpha\in G} (-b_\alpha)\neq 0.\]
Existence of algebraically free sequences of length $\alpha$ is equivalent to
the existence of free sequences of length $\alpha$ in the space ultrafilters
$\Ult(\bB)$. 

\begin{definition}
\label{6.1}
1)\quad {\em A good parameter} is a tuple $S=(\mu,\lambda,\bar{\chi})$ such that
$\mu,\lambda$ are cardinals satisfying 
\[\mu=\mu^{<\mu}<\cf(\lambda)<\lambda\quad\mbox{and}\quad (\forall\alpha<\cf(
\lambda))(\forall\xi<\mu)(\alpha^\xi<\cf(\lambda))\]
and $\bar{\chi}=\langle\chi_i: i<\cf(\lambda)\rangle$ is a strictly increasing
sequence of regular cardinals such that $\cf(\lambda)<\chi_0$, $(\forall
i<\cf(\lambda))(\chi_i^{<\mu}=\chi_i)$ and $\lambda=\sup\limits_{i<\cf(
\lambda)}\chi_i$.
\smallskip

\noindent 2)\quad Let $S=(\mu,\lambda,\bar{\chi})$ be a good parameter. Put
$\xs=\{(i,\xi): i<\cf(\lambda)\ \&\ 0\leq\xi\leq\chi_i^+\}$ and define a
forcing notion $\qs$ as follows. 
\medskip

\noindent{\bf A condition} is a tuple $p=\langle \gamma^p,w^p,u^p,\langle
f^p_{i,\xi,\alpha}: (i,\xi)\in u^p,\ \alpha<\gamma^p\rangle\rangle$ such that
\begin{enumerate}
\item[(a)] $\gamma^p<\mu$,\quad $w^p\in [\cf(\lambda)]^{\textstyle<\!\mu}$,
\quad $u^p\in [\xs]^{\textstyle<\!\mu}$, 
\item[(b)] $(\forall i\in w^p)((i,0),(i,\chi_i^+)\in u^p)$ and if $(i,\xi)\in
u^p$ then $i\in w^p$,
\item[(c)] for $(i,\xi)\in u^p$ and $\alpha<\gamma^p$,\qquad $f^p_{i,\xi,
\alpha}:u^p\longrightarrow 2$ is a function such that

if $\zeta<\xi$ then $f^p_{i,\xi,\alpha}(i,\zeta)=0$,\quad if $\xi\leq\zeta\leq
\chi_i^+$ then $f^p_{i,\xi,\alpha}(i,\zeta)=1$,\quad  and
\[f^p_{i,\xi,\alpha}\rest (u^p\setminus\{i\}\times\chi_i^+)=f^p_{i,0,\alpha}
\rest (u^p\setminus\{i\}\times\chi_i^+);\]
\end{enumerate}
{\bf the order}\quad is given by\quad $p\leq q$ if and only if $\gamma^p\leq
\gamma^q$, $w^p\subseteq w^q$, $u^p\subseteq u^q$, $f^p_{i,\xi,\alpha}
\subseteq f^q_{i,\xi,\alpha}$ (for $(i,\xi)\in u^p$, $\alpha<\gamma^p$) and 
\[(\forall (i,\xi,\alpha)\in u^q\times\gamma^q)(f^q_{i,\xi,\alpha}\rest u^p\in
\{f^p_{j,\zeta,\beta}:(j,\zeta,\beta)\in u^p\times\gamma^p\}\cup\{{\bf 0}_{
u^p}\}).\]

\noindent 3)\quad We say that conditions $p,q\in\qs$ are isomorphic if
$\gamma^p=\gamma^q$, $\otp(w^p)=\otp(w^q)$ and there is a bijection $H:u^p
\longrightarrow u^q$ (called {\em the isomorphism from $p$ to $q$}) such that
if $H_0:w^p\longrightarrow w^q$ is the order preserving mapping then:
\begin{enumerate}
\item[$(\alpha)$] $H(i,\xi)=(H_0(i),\zeta)$ for some $\zeta$,
\item[$(\beta)$]  for each $i\in w^p$, the mapping 
\[H^i:\{\xi\leq\chi_i^+:(i,\xi)\in u^p\}\longrightarrow\{\zeta\leq\chi_{H_0(
i)}^+:(H_0(i),\zeta)\in u^q\}\]
given by $H(i,\xi)=(H_0(i),H^i(\xi))$ is the order preserving isomorphism,
\item[$(\gamma)$] $(\forall\alpha<\gamma^p)(\forall (i,\xi)\in u^p)(f^p_{i,
\xi,\alpha}=f^q_{H(i,\xi),\alpha}\comp H)$. 
\end{enumerate}
\end{definition}

\noindent {\bf Remark:}\quad Variants of the forcing notion $\qs$ are used in
\cite{RoSh:651} to deal with attainment problems for equivalent definitions of
$\hd$, $\hL$. 

\begin{proposition}
\label{6.2}
Let $S=(\mu,\lambda,\bar{\chi})$ be a good parameter. Then $\qs$ is a
$\mu$--complete $\mu^+$--cc forcing notion.
\end{proposition}

\Proof Easily $\qs$ is $\mu$--closed. To show the chain condition suppose that
$\cA\subseteq \qs$ is of size $\mu^+$. Since $\mu^{<\mu}=\mu$ we may apply
standard cleaning procedure and find isomorphic conditions $p,q\in\cA$ such
that if $H:u^p\longrightarrow u^q$ is the isomorphism from $p$ to $q$ and
$H_0: w^p\longrightarrow w^q$ is the order preserving mapping then 
\begin{itemize}
\item $H_0\rest w^p\cap w^q$ is the identity on $w^p\cap w^q$, and
\item $H\rest u^p\cap u^q$ is the identity on $u^p\cap u^q$.
\end{itemize}
Next put $\gamma^r=\gamma^p=\gamma^q$, $w^r=w^p\cup w^q$, $u^r=u^p\cup u^q$. 
For $(i,\xi)\in u^r$ and $\alpha<\gamma^r$ we define $f^r_{i,\xi,\alpha}$ as
follows: 
\begin{enumerate}
\item if $(i,\xi)\in u^p$, $i\in w^p\setminus w^q$ then $f^r_{i,\xi,\alpha}=
f^p_{i,\xi,\alpha}\cup f^q_{H_0(i),0,\alpha}$,
\item if $(i,\xi)\in u^q$, $i\in w^q\setminus w^p$ then $f^r_{i,\xi,\alpha}=
f^p_{{H_0}^{-1}(i),0,\alpha}\cup f^q_{i,\xi,\alpha}$, 
\item if $i\in w^p\cap w^q$ then 
\[f^r_{i,\xi,\alpha}=(f^p_{i,0,\alpha}\cup f^q_{i,0,\alpha})\rest (u^r
\setminus\{i\}\times\chi_i^+)\cup {\bf 0}_{(\{i\}\times [0,\xi))\cap u^r}\cup
{\bf 1}_{(\{i\}\times [\xi,\chi_i^+])\cap u^r}.\] 
\end{enumerate}
Checking that $r\stackrel{\rm def}{=}\langle \gamma^r,w^r,u^r,\langle f^r_{i,
\xi}: (i,\xi)\in u^r\rangle\rangle\in\qs$ is a condition stronger than both
$p$ and $q$ is straightforward. \QED
\medskip

For a condition $p\in\qs$ let $\bB_p$ be the Boolean algebra $\bB_{(u^p,F^p)}$
for $F^p=\{f^p_{i,\xi,\alpha}:(i,\xi)\in u^p,\ \alpha<\gamma^p\}\cup\{{\bf
0}_{u^p}\}$ (see \ref{1.1}). Naturally we define  $\qs$--names $\dot{\bB}_S^*$
and $\dot{f}_{i,\xi,\alpha}$ (for $i<\cf(\lambda)$, $\xi<\chi_i^+$, $\alpha<
\mu$) by:
\[\forces_{\qs}\mbox{`` }\dot{\bB}_S^*=\bigcup\{\bB_p\!: p\in\Gamma_S\},\quad
\dot{f}_{i,\xi,\alpha}=\bigcup\{f^p_{i,\xi,\alpha}\!: (i,\xi,\alpha)\in
u^p\times\gamma^p,\ p\in\Gamma_{\qs}\}\mbox{ ''.}\]
Further, let $\dot{\bB}_S$ be the $\qs$--name for the subalgebra $\langle
x_{i,\xi}:i<\cf(\lambda), \xi<\chi_i^+\rangle_{\dot{\bB}^*_S}$ of
$\dot{\bB}^*_S$. 

\begin{proposition}
\label{6.3}
Assume $S=(\mu,\lambda,\bar{\chi})$ is a good parameter. Then in $\V^{\qs}$:
\begin{enumerate}
\item $\dot{f}_{i,\xi,\alpha}:\xs\longrightarrow 2$ (for $\alpha<\mu$ $i<\cf(
\lambda)$ and $\xi\leq\chi_i^+$),
\item $\dot{\bB}_S^*$ is the Boolean algebra $\bB_{(\xs,\dot{F})}$, where
$\dot{F}=\{\dot{f}_{i,\xi,\alpha}:(i,\xi)\in\xs,\ \alpha<\mu\}$, 
\item for each $i<\cf(\lambda)$, the sequence $\langle -x_{i,\xi}: \xi<\chi_i^+
\rangle$ is (algebraically) free in the algebra $\dot{\bB}_S$,
\item ${\bf 0}_{\xs}\in\cl(\dot{F})$, so it determines a homomorphism from
$\dot{\bB}_S^*$ to 2 (so ultrafilter). Its restriction ${\bf 0}_{\xs}\rest
\dot{\bB}_S$ has tightness $\lambda$.
\end{enumerate}
\end{proposition}

\Proof 1)--3) Should be clear.

\noindent 4) First note that if $p\in\qs$ and $i<\cf(\lambda)$ then there is a
condition $q\in\qs$ stronger than $p$ and such that 
\[(\exists\alpha<\gamma^q)(f^q_{i,0,\alpha}\rest (u^p\setminus\{i\}\times(
\chi_i^++1))\equiv 0).\]
Hence we immediately conclude that ${\bf 0}_{\xs}\in\cl(\dot{F})$. Now we look
at the restriction ${\bf 0}_{\xs}\rest \dot{\bB}_S$. First fix $i<\cf(
\lambda)$ and let $\dot{Y}_i=\{\dot{f}_{i,\xi,\alpha}\rest\dot{\bB}_S:\xi<
\chi^+_i,\ \alpha<\mu\}$ (so $\dot{Y}_i$ is a family of homomorphisms from
$\dot{\bB}_S$ to 2 and it can be viewed as a family of ultrafilters on
$\dot{\bB}_S$). It follows from the previous remark (and \ref{6.1}(2c)) that
${\bf 0}_{\xs}\rest\dot{\bB}_S\in\cl(\dot{Y_i})$. We claim that ${\bf 0}_{\xs}
\rest\dot{\bB}_S$ is not in the closure of any subset of $\dot{Y}_i$ of size
less than $\chi_i^+$. So assume that $\dot{X}$ is a $\qs$--name for a subset of
$\dot{Y}_i$ such that $\forces|\dot{X}|\leq\chi_i$ (and we will think that
$\forces \dot{X}\subseteq\chi_i^+\times\mu$). Since $\qs$ satisfies the
$\mu^+$--cc we find $\xi<\chi_i^+$ such that $\forces\dot{X}\subseteq\xi\times
\mu$. Now note that \ref{6.1}(2c) implies that $\forces(\forall(\zeta,\alpha)
\in\dot{X})(\dot{f}_{i,\zeta,\alpha}(i,\xi)=1)$, so $\forces {\bf 0}_{\xs}
\rest\dot{\bB}_S\notin\cl(\dot{X})$. Hence the tightness of the ultrafilter
${\bf 0}_{\xs}\rest\dot{\bB}_S$ is $\lambda$. \QED

\begin{theorem}
\label{6.4}
Assume that $S=(\mu,\lambda,\bar{\chi})$ is a good parameter. Then in
$\V^{\qs}$:  
\begin{enumerate}
\item there is no algebraically free sequence of length $\lambda$ in
$\dot{\bB}_S$,
\item if $\dot{I}$ is an ideal in $\dot{\bB}_S$ then $\Dep(\dot{\bB}_S/
\dot{I})<\lambda$. 
\end{enumerate}
\end{theorem}

\Proof 1)\quad Assume that $\langle\dot{b}_\alpha:\alpha<\lambda\rangle$ is a
$\qs$--name for a $\lambda$--sequence of elements of $\dot{\bB}_S$ and $p\in
\qs$. For each $i<\cf(\lambda)$ and $\xi<\chi_i^+$ choose a condition $p_{i,
\xi}\in\qs$ stronger than $p$, a finite set $v_{i,\xi}\subseteq u^{p_{i,\xi}}$
and a Boolean term $\tau_{i,\xi}$ such that 
\[p_{i,\xi}\forces_{\qs}\dot{b}_{\chi_i+\xi}=\tau_{i,\xi}(x_{j,\zeta}:(j,
\zeta)\in v_{i,\xi}).\]
Let us fix $i<\cf(\lambda)$ for a moment. Applying $\Delta$--lemma arguments
and standard cleaning (and using the assumption that $\chi_i^{<\mu}=\chi_i=\cf(
\chi_i)$) we may find a set $Z_i\in [\chi_i^+]^{\textstyle \chi_i^+}$ such
that 
\begin{enumerate}
\item[$(\alpha)_i$] all conditions $p_{i,\xi}$ for $\xi\in Z_i$ are
isomorphic,
\item[$(\beta)_i$]  $\{u^{p_{i,\xi}}:\xi\in Z_i\}$ forms a $\Delta$--system
with heart $u_i$,
\item[$(\gamma)_i$] if $\xi_0,\xi\in Z_i$ and $H:u^{p_{i,\xi_0}}
\longrightarrow u^{p_{i,\xi_1}}$ is the isomorphism from $p_{i,\xi_0}$ to
$p_{i,\xi_1}$ then $H[v_{i,\xi_0}]=v_{i,\xi_1}$ and $H\rest u_i$ is the
identity on $u_i$,
\item[$(\delta)_i$] $\tau_{i,\xi}=\tau_i$ (for each $\xi\in Z_i$),
\item[$(\varepsilon)_i$] $u^{p_{i,\xi_0}}\cap\{(j,\zeta): j<i\ \&\ \zeta<
\chi_j^+\}= u^{p_{i,\xi_1}}\cap\{(j,\zeta): j<i\ \&\ \zeta<\chi_j^+\}$
whenever $\xi_0,\xi_1\in Z_i$.
\end{enumerate}
Apply the cleaning procedure and $\Delta$--lemma again to get a set $J\in
[\cf(\lambda)]^{\textstyle \cf(\lambda)}$ such that
\begin{enumerate}
\item[$(\alpha)^*$] if $i_0,i_1\in J$, $\xi_0\in Z_{i_0}$, $\xi_1\in Z_{i_1}$
then the conditions $p_{i_0,\xi_0},p_{i_1,\xi_1}$ are isomorphic,
\item[$(\beta)^*$]  $\{u_i:i\in J\}$ forms a $\Delta$--system with heart $u^*$,
\item[$(\gamma)^*$] if $i_0,i_1\in J$, $\xi_0\in Z_{i_0}$, $\xi_1\in Z_{i_1}$
and $H:u^{p_{i_0,\xi_0}}\longrightarrow u^{p_{i_1,\xi_1}}$ is the isomorphism
from $p_{i_0,\xi_0}$ to $p_{i_1,\xi_1}$ then $H[v_{i_0,\xi_0}]=v_{i_1,\xi_1}$,
$H[u_{i_0}]=u_{i_1}$ and $H\rest u^*$ is the identity on $u^*$, 
\item[$(\delta)^*$] $\tau_i=\tau$ (for $i\in J$)
\end{enumerate}
(remember the assumptions on $\cf(\lambda)$ in \ref{6.1}(1)). Now choose $i_0
\in J$ such that $\sup\{i<\cf(\lambda):(i,0)\in u^*)\}<i_0$ and pick $\xi^0_0,
\xi^0_1\in Z_{i_0}$, $\xi^0_0<\xi^0_1$. Next take $i_1\in J$ such that 
\[i_1>i_0+\sup\{i<\cf(\lambda): (i,0)\in u^{p_{i_0,\xi^0_0}}\cup
u^{p_{i_0,\xi^0_1}}\}\]
and $u_{i_1}\cap (u^{p_{i_0,\xi^0_0}}\cup u^{p_{i_0,\xi^0_1}})=u^*$. Finally
pick $\xi^1_0,\xi^1_1\in Z_{i_1}$ such that $\xi^1_0<\xi^1_1$ and, for
$\ell<2$,  
\[u^{p_{i_1,\xi^1_\ell}}\cap (u^{p_{i_0,\xi^0_0}}\cup u^{p_{i_0,\xi^0_1}})=
u^*.\]
To make our notation somewhat simpler let $p^k_\ell=p_{i_k,\xi^k_\ell}$,
$\tau^k_\ell=\tau(x_{j,\zeta}: (j,\zeta)\in v_{i_k,\xi^k_\ell})$ (for $k,\ell
<2$) and let $H^{k_0,\ell_0}_{k_1,\ell_1}:u^{p^{k_0}_{\ell_0}}\longrightarrow
u^{p^{k_1}_{\ell_1}}$ be the isomorphism from $p^{k_0}_{\ell_0}$ to
$p^{k_1}_{\ell_1}$ (for $k_0,k_1,\ell_0,\ell_1<2$). 
 
It follows from the choice of $i_k,\xi^k_\ell$ that:
\begin{enumerate}
\item[(i)]\ \ \ if $(i,0)\in u^*$, $k<2$, $\xi<\chi_i^+$ then\quad $(i,\xi)\in
u^{p^k_0}\ \Leftrightarrow\ (i,\xi)\in u^{p^k_1}$,
\item[(ii)]\ \  if $i\in (w^{p^0_0}\cup w^{p^0_1})\cap (w^{p^1_0}\cup
w^{p^1_1})$ then $(i,0)\in u^*$. 
\end{enumerate}

Now we are defining a condition $q$ stronger than all $p^k_\ell$. So we put
$\gamma^q=\gamma^{p^0_0}$, $w^q=w^{p^0_0}\cup w^{p^0_1}\cup w^{p^1_0}\cup
w^{p^1_1}$, $u^q=u^{p^0_0}\cup u^{p^0_1}\cup u^{p^1_0}\cup u^{p^1_1}$, and for
$(j,\zeta)\in u^q$ and $\alpha<\gamma^q$ we define $f^q_{j,\zeta,\alpha}:u^q
\longrightarrow 2$ in the following manner. We declare that 
\[f^q_{j,\zeta,\alpha}\rest(\{j\}\times [0,\zeta))\cap u^q\equiv 0\quad\mbox{
and }\quad f^q_{j,\zeta,\alpha}\rest (\{j\}\times [\zeta,\chi_j^+] )\cap u^q
\equiv 1,\]
and now we define $f^q_{j,\zeta,\alpha}$ on $u^q\setminus (\{j\}\times [0,
\chi_j^+])$ letting: 
\medskip

\noindent{--} if $(j,0)\in u^*$ then 
\[f^q_{j,\zeta,\alpha}\supseteq\bigcup_{\ell,k<2} f^{p^k_\ell}_{j,0,\alpha}
\rest(u^{p^k_\ell}\setminus\{j\}\times \chi_j^+),\]
[note that in this case we have: $f^q_{j,\zeta,\alpha}(\tau^k_0)=f^q_{j,\zeta,
\alpha}(\tau^k_1)$ for $k=0,1$]

\noindent{--} if $(j,0)\in u_{i_0}\setminus u^*$ then 
\[f^q_{j,\zeta,\alpha}\supseteq\bigcup_{\ell<2} f^{p^0_\ell}_{j,0,\alpha}\rest
(u^{p^0_\ell}\setminus\{j\}\times\chi_j^+)\cup\bigcup_{\ell<2} f^{p^1_\ell}_{
H^{0,0}_{1,0}(j,0),\alpha},\]
[note that then $f^q_{j,\zeta,\alpha}(\tau^1_0)=f^q_{j,\zeta,\alpha}(\tau^1_1
)$]

\noindent{--} if $(j,0)\in u^{p^k_\ell}\setminus \bigcup\{u^{p^{k'}_{\ell'}}:
(k',\ell')\neq (k,\ell),\ k',\ell'<2\}$ then 
\[f^q_{j,\zeta,\alpha}=\bigcup_{k',\ell'<2} f^{p^{k'}_{\ell'}}_{H^{k,\ell}_{
k',\ell'}(j,\zeta),\alpha},\]
[again, $f^q_{j,\zeta,\alpha}(\tau^1_0)=f^q_{j,\zeta,\alpha}(\tau^1_1)$]

\noindent{--} if $(j,0)\in u_{i_1}\setminus u^*$ and, say, $(j,\zeta)\in u^{
p^1_0}$ then let $j^*\in w^{p^0_0}$ be the isomorphic image of $j$ (in the
isomorphism from $p^1_0$ to $p^0_0$). Choose $\zeta^*<\chi^+_{j^*}$ such that,
if possible then, $f^{p^0_0}_{j^*,\zeta^*,\alpha}(\tau^0_0)=0$ (if there is no
such $\zeta^*$ take $\zeta^*=0$). Let $\zeta'=\min\{\xi:(j^*,\xi)\in u^{p^0_1}
\ \&\ \xi\geq\zeta^*\}$ and
\[f^q_{j,\zeta,\alpha}\supseteq f^{p^0_0}_{j^*,\zeta^*,\alpha}\cup f^{p^1_0}_{
j^*,\zeta',\alpha}\cup\bigcup_{\ell<2} f^{p^1_\ell}_{j,0,\alpha}\rest
(u^{p^1_1}\setminus\{j\}\times\chi_j^+)\]
[note that $f^q_{j,\zeta,\alpha}(\tau^0_0)\leq f^q_{j,\zeta,\alpha}(\tau^1_1)$].
\medskip

\noindent It is a routine to check that $q=\langle\gamma^q,w^q,u^q,\langle
f^q_{j,\zeta,\alpha}: (j,\zeta)\in u^q,\ \alpha<\gamma^q\rangle\rangle\in\qs$
is a condition stronger than all $p^k_\ell$. It follows from the remarks on
$f^q_{j,\zeta,\alpha}(\tau^1_1)$ we made when we defined
$f^q_{j,\zeta,\alpha}$ that, by \ref{1.2}, $\bB_q\models\tau^0_0\wedge\tau^0_1
\wedge\tau^1_0\leq\tau^1_1$. Hence we conclude that $q$ forces that the
sequence $\langle\dot{b}_\alpha:\alpha<\lambda\rangle$ is not free as
witnessed by $\{\chi_{i_0}+\xi^0_0,\chi_{i_0}+\xi^0_1,\chi_{i_1}+\xi^1_0\}$
and $\{\chi_{i_1}+\xi^1_1\}$.
\medskip

\noindent 2)\quad Suppose that $\dot{I}$ is a $\qs$--name for an ideal in
$\dot{\bB}_S$ and $p\in\qs$ is such that $p\forces_{\qs}$`` $\Dep(\dot{\bB}_S/
\dot{I})=\lambda$ ''. Then for each $i<\cf(\lambda)$ we find a $\qs$--name
$\langle\dot{b}_{i,\xi}: \xi<\chi_i^+\rangle$ for a sequence of elements of
$\dot{\bB}_S$ such that 
\[q\forces_{\qs}\mbox{`` }(\forall\xi<\zeta<\chi^+_i)(0/\dot{I}<\dot{b}_{i,
\xi}/\dot{I}<\dot{b}_{i,\zeta}/\dot{I})\mbox{ ''.}\]
Repeat the procedure applied in the previous clause, now with $\dot{b}_{i,
\xi}$ instead of $\dot{b}_{\chi_i+\xi}$ there, and get $i_0,i_1,\xi^0_0,
\xi^0_1,\xi^1_0,\xi^1_1$ as there (and we use the same notation $p^k_\ell,
\tau^k_\ell, H^{k_0,\ell_0}_{k_1,\ell_1}$ as before). Now we define a
condition $q$ stronger than all the $p^k_\ell$. Naturally we let $\gamma^q=
\gamma^{p^0_0}$, $w^q=w^{p^0_0}\cup w^{p^0_1}\cup w^{p^1_0}\cup w^{p^1_1}$, 
$u^q=u^{p^0_0}\cup u^{p^0_1}\cup u^{p^1_0}\cup u^{p^1_1}$. Suppose $(j,\zeta)
\in u^q$ and $\alpha<\gamma^q$. We define $f^q_{j,\zeta,\alpha}:u^q
\longrightarrow 2$ declaring that 
\[f^q_{j,\zeta,\alpha}\rest(\{j\}\times [0,\zeta))\cap u^q\equiv 0\quad\mbox{
and }\quad f^q_{j,\zeta,\alpha}\rest (\{j\}\times [\zeta,\chi_j^+] )\cap u^q
\equiv 1,\]
and: 
\medskip

\noindent{--} if $(j,0)\in u^*$ then $f^q_{j,\zeta,\alpha}\supseteq
\bigcup\limits_{\ell,k<2} f^{p^k_\ell}_{j,0,\alpha}\rest(u^{p^k_\ell}\setminus
\{j\}\times \chi_j^+)$,

\noindent{--} if $(j,0)\in u^{p^k_\ell}$ but $(j,0)\notin u^{p^{k'}_{\ell'}}$
for $(k',\ell')\neq (k,\ell)$ then 
\[f^q_{j,\zeta,\alpha}=\bigcup\limits_{k',\ell'<2} f^{p^{k'}_{\ell'}}_{H^{k,
\ell}_{k',\ell'}(j,\zeta),\alpha},\]

\noindent{--} if $(j,0)\in u_{i_1}\setminus u^*$ then 
\[f^q_{j,\zeta,\alpha}\supseteq\bigcup_{\ell<2} f^{p^1_\ell}_{j,0,\alpha}\rest
(u^{p^1_\ell}\setminus\{j\}\times\chi_j^+)\cup\bigcup_{\ell<2} f^{p^0_\ell}_{
H^{1,0}_{0,0}(j,0),\alpha},\]

\noindent{--} if $(j,0)\in u_{i_0}\setminus u^*$ then first take $\xi^\ell=
\min\{\xi\leq\chi_j^+: (j,\xi)\in u^{p^0_\ell}\ \&\ \zeta\leq\xi\}$ (for $\ell
<2$) and next put 
\[f^q_{j,\zeta,\alpha}=f^{p^0_0}_{j,\xi^0,\alpha}\cup f^{p^0_1}_{j,\xi^1,
\alpha}\cup f^{p^1_0}_{H^{0,1}_{1,0}(j,\xi^1),\alpha}\cup f^{p^1_1}_{H^{0,0}_{
1,1}(j,\xi^0),\alpha}\]
[remember that $H^{0,1}_{1,0}[u_{i_0}]=H^{0,0}_{1,1}[u_{i_0}]=u_{i_1}$ and
both isomorphisms are the identity on $u^*$].
\medskip

\noindent It should be a routine to verify that $q=\langle\gamma^q,w^q,u^q,
\langle f^q_{j,\zeta,\alpha}: (j,\zeta)\in u^q,\ \alpha<\gamma^q\rangle\rangle
\in\qs$ is a condition stronger than all $p^k_\ell$. Note that the only case
when $f^q_{j,\zeta,\alpha}(\tau^0_0)\neq f^q_{j,\zeta,\alpha}(\tau^0_1)$  is
$(j,0)\in u_{i_0}\setminus u^*$. But then $f^q_{j,\zeta,\alpha}(\tau^1_0)=
f^q_{j,\zeta,\alpha}(\tau^0_1)$ and $f^q_{j,\zeta,\alpha}(\tau^1_1)=f^q_{j,
\zeta,\alpha}(\tau^0_0)$. Hence (by \ref{1.2}) $\bB_q\models\tau^0_1\wedge
(-\tau^0_0)\leq\tau^1_0\wedge(-\tau^1_1)$ and therefore $q\forces$``
$\dot{b}_{i_0,\xi^0_1}\wedge (-\dot{b}_{i_0,\xi^0_0})\leq\dot{b}_{i_1,
\xi^1_0}\wedge (-\dot{b}_{i_1,\xi^1_1})$ ''. Now, $q\forces$`` $\dot{b}_{i_1,
\xi^1_0}/\dot{I}\leq\dot{b}_{i_1,\xi^1_1}/\dot{I}$ '' so we conclude $q
\forces$`` $\dot{b}_{i_0,\xi^0_1}\wedge (-\dot{b}_{i_0,\xi^0_0})\in\dot{I}$
''. But the last statement contradicts $q\forces$`` $\dot{b}_{i_0,\xi^0_0}/
\dot{I}<\dot{b}_{i_0,\xi^0_1}/\dot{I}$ '', finishing the proof. \QED

\begin{conclusion}
It is consistent that there is a Boolean algebra $\bB$ of size $\lambda$ such
that there is an ultrafilter $x\in\Ult(\bB)$ of tightness $\lambda$, there is no
free $\lambda$--sequence in $\bB$ and $\tig(\bB)=\lambda\notin\Dep_{\rm Hs}(
\bB)$ (i.e.~no homomorphic image of $\bB$ has depth $\lambda$). \QED
\end{conclusion}

Let us note that in the universe $\V^{\qs}$ we have $2^{\cf(\lambda)}\geq
\lambda$. This is a real limitation -- we can prove that $2^{\cf(\lambda)}$
cannot be small in this context. In the proof we will use the following
theorem cited here from \cite{Sh:233}.

\begin{theorem}[see {\cite[Lemma 5.1(3)]{Sh:233}}]
\label{from233}
Assume that $\lambda=\sup\limits_{i<\cf(\lambda)}\chi_i$, $\cf(\lambda)<
\chi_i<\lambda$, $\mu=(2^{\cf(\lambda)})^+$. Let $X$ be a $T_{3\frac{1}{2}}$
topological space with a basis ${\mathcal B}$. Suppose that $\varphi$ is a
function assigning cardinal numbers to subsets of $X$ such that:
\begin{enumerate}
\item[(i)]\ \ \  $\varphi(A)\leq\varphi(A\cup B)\leq \varphi(A)+\varphi(B)+
\omega_0$\quad for $A,B\subseteq X$,
\item[(ii)]\ \   for each $i<\cf(\lambda)$ there is a sequence $\langle
u_\alpha:\alpha<\mu\rangle\subseteq {\mathcal B}$ such that
\[(\forall g:\mu\longrightarrow 2^{\cf(\lambda)})(\exists \alpha\neq\beta)(
g(\alpha)=g(\beta)\ \&\ \varphi(u_\alpha\setminus\cl_X(u_\beta))\geq
\chi_i),\]
\item[(iii)]\    for sufficiently large $\chi<\lambda$, if $\langle A_\alpha:
\alpha<\mu\rangle$ is a sequence of subsets of $X$ such that $\varphi(A_\alpha
)\leq\chi$ then $\varphi(\bigcup\limits_{\alpha<\mu}A_\alpha)\leq\chi$.
\end{enumerate}
Then there is a sequence $\langle u_i:i<\cf(\lambda)\rangle\subseteq{\mathcal B}$
such that 
\[(\forall i<\cf(\lambda))(\varphi(u_i\setminus\bigcup_{j\neq i}u_j)\geq
\chi_i).\qquad\qquad\qquad\QED\]
\end{theorem}

\begin{theorem}
\label{6.7}
Suppose that $\bB$ is a Boolean algebra satisfying $2^{\cf(\tig(\bB))}<\tig(
\bB)$. Then for some ideal $I$ on $\bB$ we have $\Dep(\bB/I)=\tig(\bB)$.
\end{theorem}

\Proof Let $\lambda=\tig(\bB)$ and let $\langle\chi_i:i<\cf(\lambda)\rangle$
be an increasing cofinal in $\lambda$ sequence of successor cardinals,
$\chi_0>\cf(\lambda)$ and let $\mu=(2^{\cf(\lambda)})^+$. Further, let $X$ be
the Stone space $\Ult(\bB)$ and thus we may think that $\bB={\mathcal B}$ is a
basis of the topology of $X$. Now define a function $\varphi$ on subsets of
$X$ by 
\[\begin{array}{lr}
\varphi(Y)=\sup\{\kappa\!:&\mbox{there are sequences }\langle y_\zeta\!:\zeta<
\kappa\rangle\subseteq Y \mbox{ and }\langle u_\zeta\!:\zeta<\kappa\rangle
\subseteq {\mathcal B}\ \\
&\mbox{such that }(\forall\zeta,\xi<\kappa)(y_\zeta\in u_\xi\ \Leftrightarrow\
\xi<\zeta)\}.
  \end{array}\] 
We are going to apply \ref{from233} to these objects and for this we should
check the assumptions there. The only not immediate demands might be (ii) and
(iii). So suppose $i<\cf(\lambda)$. Since $\chi_i<\lambda=\tig(\bB)$ we can
find a free sequence $\langle u^*_\xi:\xi<\chi_i^+\rangle\subseteq\bB$. Next,
for each $\xi<\chi_i^+$ we may choose an ultrafilter $y_\xi\in X$ such that
$(\forall\zeta<\chi_i^+)(y_\xi\in u^*_\zeta\ \Leftrightarrow\ \zeta<\xi)$. 
Now, for $\alpha<\mu$, let $u_\alpha=u^*_{\chi_i\cdot\alpha}$. Suppose $g:\mu
\longrightarrow 2^{\cf(\lambda)}$ and take any $\alpha<\beta<\mu$ such that
$g(\alpha)=g(\beta)$. Note that 
\[u_\alpha\setminus\cl_X(u_\beta)=u^*_{\chi_i\cdot\alpha}\setminus u^*_{\chi_i
\cdot\beta}\supseteq\{y_\xi:\chi_i\cdot\alpha<\xi<\chi_i\cdot(\alpha+1)\}\]
and easily $\varphi(\{y_\xi:\chi_i\cdot\alpha<\xi<\chi_i\cdot(\alpha+1)\})=
\chi_i$. Thus $\varphi(u_\alpha\setminus\cl_X(u_\beta))\geq\chi_i$ and the
demand \ref{from233}(ii) is verified. Assume now that $\mu<\chi<\lambda$ and
$A_\alpha\subseteq X$ (for $\alpha<\mu$) are such that $\varphi(
\bigcup\limits_{\alpha<\mu}A_\alpha)>\chi$. Let sequences $\langle y_\xi:\xi<
\chi^+\rangle \subseteq\bigcup\limits_{\alpha<\mu}A_\alpha$ and $\langle
u_\xi:\xi<\chi^+\rangle\subseteq\bB$ witness this. Then for some $C\in
[\chi^+]^{\textstyle\chi^+}$ and $\alpha<\mu$ we have $\langle y_\xi:\xi\in
C\rangle\subseteq A_\alpha$ and therefore $\langle y_\xi,u_\xi:\xi\in C
\rangle$ witness $\varphi(A_\alpha)\geq\chi^+$. This finishes checking the
demand \ref{from233}(iii).

So we may use \ref{from233} and we get a sequence $\langle u_i:i<\cf(\lambda)
\rangle\subseteq\bB$ such that 
\[(\forall i<\cf(\lambda))(\varphi(u_i\setminus\bigcup_{j\neq i}u_j)\geq
\chi_i).\] 
Then for each $i<\cf(\lambda)$ we may choose sequences $\langle y^i_\xi:\xi<
\chi_i\rangle\subseteq u_i\setminus\bigcup\limits_{j\neq i}u_j$ and $\langle
w^i_\xi:\xi<\chi_i\rangle\subseteq {\mathcal B}$ such that 
\[y^i_\xi\in w^i_\zeta\quad\Leftrightarrow\quad \zeta<\xi,\]
and we may additionally demand that $w^i_\xi\subseteq u_i$ (for each $\xi<
\chi_i$). Now let
\[I\stackrel{\rm def}{=}\{b\in\bB: (\forall i<\cf(\lambda))(\forall\xi<\chi_i)
(y^i_\xi\notin b)\}.\]
It should be clear that $I$ is an ideal in the Boolean algebra $\bB$
(identified with the algebra of clopen subsets of $X$). Fix $i<\cf(\lambda)$
and suppose that $\zeta<\xi<\chi_i$. By the choices of the $w^i_\xi$'s we have
$y^i_\xi\in w^i_\zeta\setminus w^i_\xi$ and no $y^i_\rho$ belongs to $w^i_\xi
\setminus w^i_\zeta$. As $w^i_\xi\subseteq u_i$ we conclude $\bB/I\models
w^i_\xi/I<w^i_\zeta/I$. Thus the sequence $\langle w^i_\xi/I: \xi<\chi_i
\rangle$ (for $i<\cf(\lambda)$) is strictly decreasing in $\bB/I$ and
consequently $\Dep(\bB/I)\geq \lambda$. Since there is $\lambda$ many
$y^i_\xi$'s only, we may easily check that there are no decreasing
$\lambda^+$--sequences in $\bB/I$ (remember the definition of $I$), finishing
the proof. \QED

\end{document}